\documentclass[a4paper,reqno,twoside]{amsart}
\usepackage{amsmath,amsfonts,amssymb,enumerate,stmaryrd,verbatim,epic,url}
\usepackage[english]{babel}
\usepackage{graphicx}

\theoremstyle{plain}
 \newtheorem{theorem}{Theorem}[section]
 \newtheorem{tthr}[theorem]{Theorem}
 \newtheorem{corr}[theorem]{Corollary}
 \newtheorem{lmr}[theorem]{Lemma}
 \newtheorem{prr}[theorem]{Proposition}
 \newtheorem{hyp}[theorem]{Conjecture}
\theoremstyle{definition}
 \newtheorem{defnr}[theorem]{Definition}
 \newtheorem{remnr}[theorem]{Remark}

\def\brm{\begin{remnr}}
\def\erm{\end{remnr}}
\def\bdr{\begin{defnr}}
\def\edr{\end{defnr}}
\def\bp{\begin{proof}}
\def\ep{\end{proof}}
\def\btr{\begin{tthr}}
\def\etr{\end{tthr}}
\def\bpr{\begin{prr}}
\def\epr{\end{prr}}
\def\bcr{\begin{corr}}
\def\ecr{\end{corr}}
\def\blr{\begin{lmr}}
\def\elr{\end{lmr}}
\def\beq{\begin{equation}}
\def\eeq{\end{equation}}
\def\bcs{\begin{cases}}
\def\ecs{\end{cases}}
\def\bhr{\begin{hyp}}
\def\ehr{\end{hyp}}

\def\G{\Gamma}
\def\lab{\label}
\def\en{\begin{enumerate}}
\def\ene{\end{enumerate}}

\def\>{\geqslant}
\def\<{\leqslant}
\def\lb{\lfloor}
\def\rb{\rfloor}

\def\quat#1{\textquotedblleft #1\textquotedblright}

\title{NOTES ON FIBONACCI PARTITIONS}

\author{F.~V.~WEINSTEIN}

\address
{Giacomettistrasse 33A, CH-3006, Switzerland.}

\email{felix.weinstein46@gmail.com}

\dedicatory{To my children Sergej, Jelena, and Maria with love}

\begin{document}

\begin{abstract}
Let $f_1=1,f_2=2$ and $f_i=f_{i-1}+f_{i-2}$ for $i>2$ be the
sequence of Fibonacci numbers.
Let $\Phi_h(n)$ be the quantity of partitions
of natural number $n$ into $h$ different Fibonacci numbers.
In terms of Zeckendorf partitions I deduce a formula for the function
$\Phi(n;t):=\sum_{h\geq 1}\Phi_h(n)t^h$, and use it to analyze
the functions $F(n):=\Phi(n;1)$ and $\chi(n):=\Phi(n;-1)$.
I obtain the least upper bound
for $F(n)$ when $f_i\<n<f_{i+1}$. In particular, it implies that
$F(n)\<\sqrt{n+1}$ for any natural $n$.

I prove also that $|\chi(n)|\<1$, and
$\mathop{\lim}\limits_{N\to\infty}\frac{1}{N}
\left(\chi^2(1)+\chi^2(2)+{\dots}+\chi^2(N)\right)=0$.

For any $k\>2$, I define a special finite set $\mathbb{G}(k)$ of
solutions of the equation $F(n)=k$; all solutions can be easily
obtained from $\mathbb{G}(k)$. This construction uses a nonstandard
representation of rational numbers as certain continued fractions and
provides with a canonical identification
$\coprod_{k\>2}\mathbb{G}(k)=\G_+$, where $\G_+$ is the monoid freely
generated by the positive rational numbers $<1$.

Let $\Psi(k)$ be the cardinality of $\mathbb{G}(k)$. I prove that,
for $i\>2k$ and $k\>2$, the interval $[f_i,f_{i+1}-2]$ contains
exactly $2\Psi(k)$ solutions of the equation $F(n)=k$ and offer a
formula for the Dirichlet generating function of the sequence
$\Psi(k)$. In addition, I study the set of minimal solutions of
the equation $F(n)=k$ as $k$ varies and I offer a
conjecture on the distribution of such solutions.
\end{abstract}

\maketitle

\markboth{F.~V.~WEINSTEIN}{NOTES ON FIBONACCI PARTITIONS}

\section{Introduction}

Most of the results of this
article appeared initially as conjectures, based on computer experiments.
Several conjectures I was unable to prove are formulated in what follows.

Let $f_0=f_1=1$ and $f_i=f_{i-1}+f_{i-2}$ for $i\>2$ be the sequence
of Fibonacci numbers\footnote{The \quat{conventional} indexing of Fibonacci numbers,
see \url{http://en.wikipedia.org/wiki/Fibonacci_number},
is different from indexing in this article.
The formulations of main statements are more neat with our indexing
than with \quat{conventional} one, cf.\cite{RN}, where the same problem
is encountered.}. A \textit{Fibonacci partition} of a natural number
$n$ is a representation $n=f_{i_1}+f_{i_2}+{\dots}+f_{i_h}$, where
$1\<i_1<i_2<{\dots}<i_h$. The numbers $f_{i_1},f_{i_2},{\dots},f_{i_h}$
are referred to as the \textit{parts} of the Fibonacci partition.

Let $\Phi_h(n)$ be the quantity of Fibonacci partitions of the natural
number $n$ with $h$ parts. Define
\[
\Phi(n;t):=\sum_{h\geq 1}\Phi_h(n)t^h,
\]
and $\Phi(0;t):=1$.
The main Theorem \ref{th_F} offers an explicit formula for
$\Phi(n;t)$ in terms of the Zeckendorf partition of $n$, see
Def.\ref{zeck}. This formula is the foundation for all subsequent
results of the article.

Furthermore, I study properties of the function
$F(n):=\Phi(n;1)$ which counts the
quantity of all Fibonacci partitions of $n$, and the function
$\chi(n):=\Phi(n;-1)$ which counts the difference between the
quantities of Fibonacci partitions of $n$ with even and odd number
of parts.
Theorem \ref{th_F} easily implies that
\[
\prod_{i=1}^\infty (1-x^{f_i})=1+\sum^\infty_{n=1}\chi(n)x^n,
\]
where $\chi(n)=0,\pm 1$.

Let $\mathbb{N}$ be the set of natural numbers (i.e., the positive integers),
and let $\mathbb{N}(k)\subset\mathbb{N}$ be the set of all solutions of the
equation $F(n)=k$. Theorem implies that $F(n)=1$ whenever $n=f_i-1$, where $i\geqslant 1$.
The main goal of Section \ref{Gen} is to offer an algorithm, which calculates
the set $\mathbb{N}(k)$, where $k\geqslant 2$.

Let me outline its construction. The main idea originated from Theorem \ref{th_F},
which may be interpreted by means of the representations of rational numbers as
\quat{nonstandard} continued functions (observe minus sign in formula \eqref{chain}).
On $\mathbb{N}$, I explicitly define, using this connection, a canonical action of the semigroup
$H=\mathbb{Z}_{\geqslant 0}\times\mathbb{Z}_2\times\mathbb{Z}_2$, where $\mathbb{Z}_{\geqslant 0}$
is the semigroup of nonnegative integers and $\mathbb{Z}_2$ is the group with two elements.

The function $F$ is invariant with respect to the $H$-action.
Observe also that $H$ freely acts on the set $\mathbb{N}\setminus\mathbb{N}(1)$.
I canonically define a fundamental domain $\mathbb{G}\subset\mathbb{N}\setminus\mathbb{N}(1)$
of this action. The elements of $\mathbb{G}$ are called the \emph{generating numbers}. Set
\[
\mathbb{G}(k):=\{n\in\mathbb{G} \mid F(n) = k\}.
\]
Then $\mathbb{G}=\coprod_{k\geqslant 2}\mathbb{G}(k)$.
The elements of $\mathbb{G}(k)$ are said to
be the \emph{$k$-generating numbers}.
The set of generating numbers possesses several
remarkable properties; it has a natural structure of a
noncommutative monoid canonically isomorphic to the
\emph{monoid freely generated by the set of positive rational
numbers $<1$}. Denote by $\times$ the monoidal multiplication.
Then, for $n1,n2\in\mathbb{G}$, we have $F(n1\times n2) = F(n1)F(n2)$.
Thus, $G(k_1)\times G(k_2)\subset G(k_1k_2)$.
Moreover, for any $k\geqslant 2$, the set $\mathbb{G}(k)$ is finite. Set
\[
\Psi(k):=|\mathbb{G}(k)| .
\]
The algorithm presented in Section \ref{Gen} makes it possible
to build the sets $\mathbb{G}(k)$. Any $n_0\in\mathbb{G}(k)$ obtained generates
an infinite set of solutions of the equation $F(n)=k$ as the
$H$-orbit of $n_0$. For different numbers $n_0\in\mathbb{G}(k)$, the corresponding
orbits do not intersect, and their union coincides
with the set $\mathbb{N}(k)$. Therefore, since the definition of
$H$-action is explicit, we can obtain explicit formulas for all
solutions of the equation $F(n)=k$ as soon as the set of $k$-generating
numbers and their corresponding Zeckendorf
partitions are known.

For example, for a prime number $k$, all solutions are
easy to find if we know the representation of any rational
number $\frac{m}{k}$,
where $m=1,2,\dots,k-1$, as a continued
fraction of the form \eqref{chain}.

It is interesting to note that in our approach a \emph{noncommutative object}
— the monoid $\mathbb{G}$ of generating numbers
with multiplication operation $\times$ quite naturally appears.
Moreover, the monoidal structure of $\mathbb{G}$ allows us to completely
describe the structure of the whole set of solutions
of the equation $F(n) = k$ as $k$ varies.

Although in this article it is not used, in Theorem \ref{G} I
give a transparent description of $\mathbb{G}$:
\[
\mathbb{G}=\left\{2\lfloor l\phi\rfloor+l\;\mid\;\text{$l\in\mathbb{N}$\,
and\, $\phi:=\frac{1}{2}\left(1+\sqrt{5}\right)$}\right\}.
\]
In this form the sequence of generating numbers was
known since at least 1992 (see \cite{SL},\,A003623).

Auxiliary Section \ref{EstD} contains some inequalities needed
in the next two sections.
In Section \ref{Stab}, first I show that the maximal $k$-generating
number is equal to $f_{2k}-2$. Using this and the $H$-action
on $\mathbb{N}$ I establish a \quat{stabilization} of quantities of Fibonacci
partitions in the following sense: if $i\>2k$, then
\[
\big|\left\{\,n\,\mid\,f_i-1\<n<f_{i+1}-1\text{~~and~~$F(n)=k$}\,\right\}\big|=
\bcs
1&\text{if $k=1$},\\
2\Psi(k)&\text{if $k>1$}.
\ecs
\]
In particular, this gives the asymptotic
\[
\big|\{\,n\leq N\,\mid\,F(n)=k\,\}\big|\;\underset{N\to\infty}\sim\;2\Psi(k)\log_\phi(N),
\qquad\text{where $\;k\>2\,$}.
\]

In Section \ref{Up}, for all $n$ such that $f_i-1\<n<f_{i+1}-1$, I establish the least upper
bound
\[
F(n)\<
\bcs
f_{\frac{i+1}{2}}&\text{\rm if\; $i\equiv 1\mod 2$},\\[1mm]
2f_{\frac{i}{2}-1}&\text{\rm if\; $i\equiv 0\mod 2$},
\ecs
\]
and list all the cases in which the equalities hold.

One corollary of this result is the inequality $F(n)\leqslant\sqrt{n+1}$
valid for any natural $n$, where the equality holds
whenever $n=f^2_i-1$.

In addition, it turns out that if $i\not\in\{1,2,3,4,6,9\}$
and $f_i-1\leqslant f_{i+1}-1$, the upper bound of $F(n)$ is
attained for exactly 2 values of $n$ if $i\equiv 1\mod 2$, and for
exactly 4 values of n if $i\equiv 0\mod 2$.
The material of Section \ref{chi} concerns the function $\chi(n)$.
Since $|\chi(n)|\leqslant 1$, the quantity of natural numbers $n$ such
that $n\leqslant N$ with $\chi(n)=\pm 1$ is equal to
\[
X(N) := \chi^2(1) + \chi^2(2)+\dots+\chi^2(N).
\]

The main result of Section \ref{chi} (Theorem \ref{th41}) establishes that
$\mathop{\lim}\limits_{N\to\infty}\frac{X(N)}{N}=0$. This is a bit
surprising since the natural numbers $n$ with $\chi(n)\neq 0$ appear rather
often. For example, $X(f_{26})=X(196418)=46299$.
The proof is based only on Theorem \ref{th_F}.
Additional information on the behaviour of sequence $\chi(n)$
is contained in Remark \ref{RemXi}.

Computing the minimal $k$-generating number, i.e., the
minimal solution of the equation $F(n) = k$, is a difficult
task. Denote this solution by $m_F(k)$. The main result of
Section \ref{Fnk} states that there exists a uniquely defined
proper subset $P\subset\mathbb{N}$ such that, for any natural $k\geqslant 2$,
there is a unique decomposition
\[
m_F(k) = m_F (k_1)\times m_F (k_2)\times\dots\times m_F (k_{r(k)}),
\]
where $k_1, k_2,\dots, k_{r(k)}\in P$ and $r(k) > 1$
whenever $k\not\in P$

The numbers -- elements of $P$ will be called $F$-primitive.
The set $P$ contains prime numbers, Fibonacci numbers,
and certain other numbers. In general, the nature of $P$
is completely obscure to me. To make it clearer, I analyzed
the numerical data presented in \cite{SL}, A013583.
As a result, I was able to formulate in Section 8 a conjecture
concerning “distribution” of the numbers $m_F (k)$
when $k\in P$. Based on this analysis, I can also conjecture
that asymptotically we have

\[
\big|\{k\leqslant N\mid\text{$k$ is $F$-primitive}\}\big|\underset{N\to\infty}\sim cN,
\qquad\text{where\; $c\approx 0.35$.}
\]

For several remarks on the graph of the function $F$, see
Section \ref{Fgr}. The main goal of Section \ref{Fgr} is to offer a conjecture
describing the boundary of the convex hull of the set
of points $\{(n,F(n))|n\in [f_i-1, f_{i+1}-1]\}$.

In Section \ref{Add}, I study the function $\Psi(k)$, and some
arithmetical functions naturally related to it. (As an aside,
observe that $\Psi(k)$ was first defined by R. Munafo in 1994
as the quantity of continental $\mu$-atoms of period $k$ in the
Mandelbrot set, see \cite{SL}, A006874.)

I show how one can recursively calculate $\Psi(k)$. Furthermore,
in Theorem \ref{dgft} I prove a formula which
implies the following expression for the Dirichlet generating
function of the sequence $\Psi(k)$ in terms of the Riemann
$\zeta$ -function:
\beq\lab{zet}
1 +\sum_{k=2}^\infty\frac{\Psi(k)}{k^s}=\left(2-\frac{\zeta(s-1)}{\zeta(s)}\right)^{-1}.
\eeq
Additional arithmetical functions considered in Section \ref{Add}
appear as generating functions for orbit sets of some
naturally defined subgroups in the automorphism group
of the monoid G. In passing, three known sequences
appear here: ordered Bell numbers, Bell numbers, and
numbers of the unordered multiplicative partitions. This
part of Section \ref{Add} does not contain results directly related
to Fibonacci partitions, but is a collection of ad hoc
remarks.

Most of the above-mentioned results were obtained at
the end of 1994 and presented in 1995 as an internal report
of Bern University (Report No. 329/FW 1, January 1995,
pp. 1–29). In 2003, I put a revised version of the report in arXiv.

Here, essential additions to the 2003-version are the
main result of Section \ref{Up}, result of Section \ref{Fnk}, and formula
\eqref{zet}. Among other things, the result of Section \ref{Up} implies
the inequality $F(n)\leqslant\sqrt{n+1}$
in a more natural and clear
way as compared with the earlier proof. I also give an
explicit algorithm for deriving all solutions of the equation
$F(n) = k$; the 2003-version contains only an implicit
form of the algorithm.

Two of the mentioned results were published earlier:
the inequality $|\chi(n)|\leqslant 1$ was obtained by N. Robbins in
\cite{RN}; Theorem \ref{th41} was obtained by F. Ardilla in
\cite{ARP}. My proofs of these claims are based on Theorem \ref{th_F}
and differ from the proofs in \cite{RN} and \cite{ARP}.

\emph{History}. Let me shortly outline some earlier publications
on Fibonacci partitions. Mainly, they appeared in
“Fibonacci Quarterly” during 1960s. Almost all of them
(e.g., \cite{HOG,KL,CAR}) deal
with the calculation of $F(n)$ for some particular values of
$n$ like $n=f_i,f_i-1,f_i^2-1, f_{2i}f_{2i+1}$, etc.
In \cite{KL}, the solutions of the equation $F(n) = k$ for $k = 1, 2, 3$
were also obtained. During that time the only general
result on $F(n)$ was presented in the article by Carlitz
\cite{CAR}, where a recurrent formula for $F(n)$ was
established by means of the Zeckendorf partition. But,
as Carlitz notes, it is too complicated for practical usage.
Later on, two other general results were obtained in the
already quoted articles \cite{RN} and \cite{ARP}.
To the best of my knowledge, the cited articles cover all
currently published essential information on Fibonacci
partitions.

{\it Acknowledgements.}\; I am grateful to N.~J.~A.~Sloane who read
my 1995-report and included some of the integer sequences from it in
\cite{SL}. I also would like to thank R.~Munafo for a kind letter,
where he explained his result, N.~Robbins for sending me a draft of
his article \cite{RN}, R.~P.~Stanley who pointed me to the paper \cite{ARP}, and
A.~M.~Vershik for useful discussions.
\smallskip

\noindent
\textbf{Notation:}
\smallskip

$\mathbb{Z}$ -- the ring of integer numbers.
\vspace{1mm}

$\mathbb{N}$ -- the commutative semigroup (by multiplication) of
natural (positive integer) numbers.
\vspace{1mm}

$\mathbb{Q}_{>0}$ -- the set of positive rational numbers.
\vspace{1mm}

$\mathbb{Q}_{(0,1)}$ -- the set of positive rational numbers $<1$.
\vspace{1mm}

$|M|$ -- the cardinality of the set $M$. \vspace{1mm}

$M_1\sqcup M_2$ -- the disjoint union of sets $M_1$ and $M_2$.
\vspace{1mm}

$\varphi(n)$ -- the quantity of naturals $r<n$ relatively prime to
$n$ (the \textit{Euler totient function}). \vspace{1mm}

$\lfloor a\rfloor$ and $\lceil a\rceil $ denote the maximal integer
$\<a$ and the minimal integer $\>a$, respectively. \vspace{1mm}

For a finite set of numbers $I$, define $I^-:=\min(I)$ and
$I^+:=\max(I)$. \vspace{1mm}

For any $a\in\mathbb{Z}$, define $ \beta(a):=
\bcs 0&\text{if $a\equiv 0\mod 2$},\\
1&\text{if $a\equiv 1\mod 2$}. \ecs $ \vspace{1mm}

$m*a$ for any $m\in\mathbb{N}$ is a shorthand for $a, \dots, a$,
where $a$ is repeated $m$ times. \vspace{1mm}

$[a,b]:=\left\{n\in\mathbb{Z}\,\mid\,a,b\in\mathbb{Z}\;\;\text{and}\;\;a\<n\<b\right\}$.

\section{\bf Quantity of Fibonacci partitions}\lab{main}

In this section we formulate and prove the main result of the paper -- Theorem \ref{th_F}.
First, let us introduce necessary definitions.
\bdr
A (strict) \textit{partition} is a finite ordered set $I=\langle i_1,i_2,\dots, i_h\rangle$
of natural numbers such that $i_1<i_2<{\dots}<i_h$.
These numbers are called \textit{parts of the partition $I$}.
The numbers
\beq\lab{len}
\|I\|=i_1+i_2+{\dots}+ i_h\qquad\text{and}\qquad |I|=h
\eeq
are called the
\textit{degree} and \textit{length} of $I$, respectively.
\edr

\bdr
The \textit{union} $I_1\sqcup I_2$ of partitions $I_1$ and
$I_2$ with $I_1^+<I_2^-$ is the partition whose set of parts is the
union of the sets of parts of $I_1$ and $I_2$.
\edr

\bdr
A 2-\textit{partition} is a partition $\langle i_1,i_2,\dots, i_h\rangle$ such
that $i_{r+1}-i_r\>2$ for all $r$ such that $1\<r<h$. The set of all
2-partitions is denoted by $\mathbb{P}$.
\edr

\bdr
A 2-partition $\langle i_1,i_2,\dots, i_h\rangle$ is called \textit{simple}
if $\beta(i_1)=\beta(i_2)={\dots}=\beta(i_h)$.

For any 2-partition $I$, there exists a unique decomposition
$I=I_1\sqcup I_2\sqcup{\dots}\sqcup I_r$, called the
\emph{canonical decomposition} of $I$, where $I_1,I_2{\dots},I_r$ are simple 2-partitions
and $\beta(I_m^+)\neq\beta(I_{m+1}^-)$ for all $m$ such that $1\<m<r$.
\edr

For example, $\langle 3,5,10,12\rangle=\langle 3,5\rangle\sqcup\langle10,12\rangle$.

\bdr
For any $2$-partition $I=\langle i_1,i_2,\dots,i_h\rangle$, define
a vector $a(I)=(a_1,a_2,\dots,a_h)$, where
\[
a_l=
\bcs
\left\lb\frac{i_1-1}{2}\right\rb +1       & \text{if $\;l=1$,}\\[1mm]
\left\lb\frac{i_l-i_{l-1}}{2}\right\rb +1 & \text{if $\;1<l\<h$}.
\ecs
\]
\edr
Note that $a_1\>1$ and $a_l\>2$ if $l\>2$.
For example, $a\left(\langle 3,5,10,12\rangle\right)=(2,2,3,2)$.



\bdr\lab{defPoly}
Let $I$ be a $2$-partition, let $I=I_1\sqcup I_2\sqcup\dots\sqcup I_r$
be the canonical decomposition, and let $a(I)=(\alpha_1(I),\alpha_2(I),\dots,\alpha_r(I))$,
where $|\alpha_m(I)|=|I_m|$ for any $m=1,2,\dots,r$.

A \emph{polyvector of $I$} is an expression of the form
\[
\alpha(I)=\alpha_1(I)\times\alpha_2(I)\times\dots\times\alpha_r(I).
\]
\edr

For example, $\alpha(\langle 3,5,10,12\rangle)=\alpha(\langle 3,5\rangle\sqcup\langle10,12\rangle)=(2,2)\times(3,2)$.

Let $I$ be a $2$-partition and let
$\alpha(I)=\alpha_1(I)\times\alpha_2(I)\times\dots\times\alpha_r(I)$
be the corresponding polyvector. Set
\[
\Delta(I;t)=\Delta(\alpha_1(I);t)\cdot\Delta(\alpha_2(I);t)\cdot{\dots}\cdot\Delta(\alpha_r(I);t),
\]
where for any vector $A=(a_1,a_2,\dots, a_m)$ with natural coordinates,
\[
\Delta(A;t)=\Delta\left(a_1,a_2,\dots, a_m;t\right):=
\begin{vmatrix}
\;\xi(a_1;t) & t^{a_2+1} & 0 & 0 & \hdots & 0 \\
\;1 & \xi(a_2;t) & t^{a_3+1} & 0 & \hdots & 0 \\
\;\vdots & \vdots & & \ddots & &\vdots \\
\;0 & 0 & \hdots & 1 &\xi(a_{m-1};t) & t^{a_m+1} \\
\;0 & 0 & \hdots & 0 & 1 &\xi(a_m;t)
\end{vmatrix}
\]
and for natural $a$ we have
\[
\xi(a;t):=t+t^2+{\dots}+ t^a.
\]
The polynomial
$\Delta(A;t)$ can be defined recursively by the formulas
\beq\lab{eqRec}
\begin{gathered}
\Delta(a_1;t):=\xi(a_1;t),\qquad\Delta(a_1,a_2;t):=\xi(a_2;t)\Delta(a_1;t)-t^{a_2+1},\\
\Delta(a_1,\dots, a_m;t):=\xi(a_m;t)
\Delta(a_1,\dots, a_{m-1};t)-t^{a_m+1}\Delta(a_1,\dots, a_{m-2};t)\qquad\text{for\; $m>2$.}
\end{gathered}
\eeq

Now let us turn to Fibonacci partitions and their connection with $2$-partitions.

\bdr
A \textit{Fibonacci partition} is a partition all whose parts
are Fibonacci numbers. For any partition $I=\langle i_1,i_2,\dots, i_h\rangle$,
define $f_I:=\langle f_{i_1},f_{i_2},\dots, f_{i_h}\rangle$.

The set of all Fibonacci
partitions of degree $n$ is denoted by $\mathbb{F}(n)$.
\edr

\bdr\lab{zeck}
A \textit{Zeckendorf partition of degree $n$} is any partition
$f_I\in\mathbb{F}(n)$, where $I$ is a 2-partition.
\edr

It is almost self-evident that
\textit{for any natural $n$, there exists a
unique Zeckendorf partition of degree $n$} (see \cite{GKP},
Sec.6.6). We will denote this partition by
\[
\text{$f_{Z(n)}$,\quad where\quad $Z(n):=\left\langle z(n,1),z(n,2),\dots, z(n,\lambda(n))\right\rangle$\;\; is a 2-partition.}
\]
Thus, $n=\|f_{Z(n)}\|$.
For example,
\[
Z(333)=\langle 3,5,10,12\rangle\qquad\text{since}\qquad 333=f_3+f_5+f_{10}+f_{12}.
\]

\bdr
A natural $n$ is called \textit{$f$-simple} if the 2-partition
$Z(n)$ is simple.
\edr

\bdr
Let $M$ be a set of partitions of equal degree, and let $m_l$ be
the quantity of these partitions of length $l$. The polynomial
$G_M(t):=\sum_l m_l t^l$ is called the \textit{generating function of $M$}.
\edr

Define
\[
\Phi(0;t)=1\qquad\text{and}\qquad
\Phi(n;t)=G_{\mathbb{F}(n)}(t)\quad\text{for integer $n>0$}.
\]
The following statement is the main result of the section:

\btr\lab{th_F}
For $n$ natural, let $\alpha(Z(n))=\alpha_1(n)\times\alpha_2(n)\times\cdots\times\alpha_r(n)$
be the polyvector of the $2$-partition $Z(n)$.
Then
\beq\lab{PhiD}
\Phi(n;t)=\Delta\left(\alpha_1(n);t\right)\cdot\Delta\left(\alpha_2(n);t\right)
\cdot{{\dots}}\cdot\Delta\left(\alpha_r(n);t\right).
\eeq
\etr

For example,
\[
\Phi(333;t)=\Delta(2,2;t)\cdot\Delta(3,2;t)=t^4+2 t^5+4 t^6+4 t^7+3 t^8+t^9.
\]

Proof of Theorem \ref{th_F} is based on several lemmas.
For brevity, until the end of this section we will write,
as a rule, $\lambda$ instead of $\lambda(n)$.

\blr\lab{L_Flem01}
\noindent
The maximal part of any partition from
$\mathbb{F}(n)$ is equal to either $f_{z(n,\lambda)}$, or $f_{z(n,\lambda)-1}$.
\elr

\bp
If $\langle f_{i_1},f_{i_2},\dots, f_{i_h}\rangle\in\mathbb{F}(n)$
and $f_{i_h}\<f_{z(n,\lambda)-2}$, then
\[
n=f_{i_1}+{\dots}+f_{i_h}\<f_1+f_2+{\dots}+f_{z(n,\lambda)-2}=f_{z(n,\lambda)}-2<n,
\]
which is impossible.
\ep

\blr\lab{L_Flem02}
Any partition from $\mathbb{F}(f_i)$ is of the form
\[
\left\langle f_{i-2s},f_{i-(2s-1)},f_{i-(2s-3)},\dots, f_{i-3}+f_{i-1}\right\rangle,
\quad\text{where}\quad s=0,1,\dots,\left\lb\frac{i-1}{2}\right\rb.
\]
In particular, $\Phi(f_i;t)=\xi\left(\left\lb\frac{i-1}{2}\right\rb+1;t\right)$.
\elr

\bp
This follows by induction on $s$ from Lemma \ref{L_Flem01}.
\ep

\blr\lab{L_Flem03}
For any partition $\left\langle
f_{i_1},f_{i_2},\dots, f_{i_h}\right\rangle\in\mathbb{F}(n)$ and
any $m$ such that $1\<m\<\lambda$, there exists an $s$ such
that $1\<s\<h$ and $\left\langle
f_{i_s},f_{i_{s+1}}{\dots},f_{i_h}\right\rangle\in\mathbb{F}\left(f_{z(n,m)}+{\dots}+f_{z(n,\lambda)}\right)$.
\elr

\bp
The induction on $\lambda$ shows that it suffices to consider
only the case where $m=\lambda$. By Lemma \ref{L_Flem01} either
$i_h=z(n,\lambda)$, or $i_h=z(n,\lambda)-1$. If $i_h=z(n,\lambda)$, then $s=h$. Let
$i_h=z(n,\lambda)-1$ and let $s\>1$ be the maximal number such that
$i_{s+1}=i_s+1$ and $i_{m+1}-i_m=2$ for all $m$ such that $s<m<h$.
Then $f_{i_s}+{\dots}+ f_{i_h}=f_{z(n,\lambda)}$.
\ep

\blr\lab{L_Flem04}
Let $z(n,1)>m\>0$. Set
\[
\mathbb{F}_m(n):=\left\{f_I\in\mathbb{F}(n)\mid I^->m\right\}.
\]
Then
\[
G_{\mathbb{F}_m(n)}(t)=\Phi\left(f_{z(n,1)-m}+f_{z(n,2)-m}+{\dots}+f_{z(n,\lambda)-m};t\right).
\]
\elr

\bp
For $\left\langle f_{i_1},\dots, f_{i_h}\right\rangle\in\mathbb{F}_m(n)$, set
$\nu_m\left(\left\langle f_{i_1},\dots, f_{i_h}\right\rangle\right)=
\left\langle f_{i_1-m},\dots, f_{i_h-m}\right\rangle$.
Let us check that
\[
\nu_m\left(\langle f_{i_1},\dots, f_{i_h}\rangle\right)\in\mathbb{F}
\left(f_{z(n,1)-m}+f_{z(n,2)-m}+{\dots}+f_{z(n,\lambda)-m}\right).
\]
By Lemma \ref{L_Flem03} there exists an $s\>1$ such that
$f_{i_s}+{\dots}+f_{i_h}=f_{z(n,\lambda)}$.
For $s=1$, the claim follows from Lemma \ref{L_Flem02}.
For $s>1$, we use the induction on $s$. The inductive hypothesis shows that
\[
f_{i_1-m}+{\dots}+ f_{i_{s-1}-m}=f_{z(n,1)-m}+{\dots}+ f_{z(n,\lambda-1)-m}
\quad\text{and}\quad
f_{i_s-m}+{\dots}+ f_{i_h-m}=f_{z(n,\lambda)-m}.
\]
This gives the claim required. Thus, the map
\[
\nu_m:\mathbb{F}_m(n)\to\mathbb{F}\left(f_{z(n,1)-m}+f_{z(n,2)-m}+{\dots}+f_{z(n,\lambda)-m}\right)
\]
is well-defined and injective. Hence, $\nu_m$ is bijective.
\ep

\blr\lab{L_Flem05}
Let $\left\langle f_{z(n,1)},{\dots},f_{z(n,\lambda_1)},f_{z(n,\lambda_1+1)},{\dots},
f_{z(n,\lambda)}\right\rangle$ and
$n_1=f_{z(n,1)}+{\dots}+f_{z(n,\lambda_1)}$, where $1\<\lambda_1<\lambda$. If
$\beta\left(z(n,\lambda_1+1)\right)\neq\beta\left(z(n,\lambda_1)\right)$, then
\beq\lab{Leit1}
\Phi(n;t)=\Phi(n_1;t)\cdot
\Phi(f_{z(n,\lambda_1+1)-z(n,\lambda_1)+1}+{\dots}+f_{z(n,\lambda)-z(n,\lambda_1)+1};t).
\eeq
\elr

\bp
For $f_I=\langle f_{i_1},f_{i_2},\dots,
f_{i_h}\rangle\in\mathbb{F}(n)$, let $f_{I_2}=\langle
f_{i_s},f_{i_{s+1}},\dots, f_{i_h}\rangle\in \mathbb{F}(n-n_1)$. Such
an $s$ exists by Lemma \ref{L_Flem03}. Then $f_{I_1}=\langle
f_{i_1},f_{i_2},\dots, f_{i_{s-1}}\rangle\in\mathbb{F}(n_1)$.

By Lemma \ref{L_Flem01} either $i_{s-1}=z(n,\lambda_1)$, or
$i_{s-1}=z(n,\lambda_1)-1$. Since $i_s>i_{s-1}$, then $i_s>z(n,\lambda_1)-1$.
Thus, any partition $f_I\in\mathbb{F}(n)$ can be uniquely written
as an element $(f_{I_1},f_{I_2})$ of the set
\[
M(n,\lambda_1):=\mathbb{F}(n_1)\times\mathbb{F}_{z(n,\lambda_1)-1}(n-n_1).
\]

This gives an imbedding $u:\mathbb{F}(n)\to M(n,\lambda_1)$.
The formula
\[
\left|(f_{I_1},f_{I_2})\right|=\left|f_{I_1}\right|+\left|f_{I_2}\right|
\]
defines a length function on
$M(n,\lambda_1)$. By Lemma \ref{L_Flem04}, we see
that the corresponding generating function
coincides with the right side of formula \eqref{Leit1}.

Let $(f_{J_1},f_{J_2})\in M(n,\lambda_1)$.
Lemma \ref{L_Flem02} implies that $\beta(i_s)=\beta(z(n,\lambda_1+1))$.
If $\beta(i_s)\neq\beta(z(n,\lambda_1))$, then $J_2^->z(n,\lambda_1)$.
Since $J_1^+\<z(n,\lambda_1)$, then $f_{J_1\sqcup J_2}\in\mathbb{F}(n)$.
Hence the map $u$ is bijective.
Therefore, the generating functions of the sets $M(n,\lambda_1)$ and $\mathbb{F}(n)$ coincide.
\ep

\blr\lab{L_Flem06}
Let $a=\frac{1}{2}\left(z(n,\lambda)-z(n,\lambda-1)\right)+1$,
where $n$ is an $f$-simple number. Then
\beq\lab{Leit2}
\Phi(n;t)=\Phi\left(n-f_{z(n,\lambda)};t\right)\cdot\Phi\left(f_{z(n,\lambda)-z(n,\lambda-1)+1};t\right)-
t^{a+1}\Phi\left(n-f_{z(n,\lambda-1)}-f_{z(n,\lambda)};t\right).
\eeq
\elr

\bp
Let
\[
M(n)=\mathbb{F}\left(n-f_{z(n,\lambda)}\right)\times\mathbb{F}_{z(n,\lambda_1)-1}\left(f_{z(n,\lambda)}\right).
\]
Define an imbedding $u:\mathbb{F}(n)\to M(n)$ as in the proof of Lemma \ref{L_Flem05}.
Applying Lemma \ref{L_Flem04} we see that the generating function of $M(n)$ is equal to
\[
\Phi\left(n-f_{z(n,\lambda)};t\right)\cdot\Phi\left(f_{z(n,\lambda)-z(n,\lambda-1)+1};t\right).
\]
There is a unique partition $f_A\in\mathbb{F}_{z(n,\lambda)-1}(f_{z(n,\lambda)})$
such that $A^-=z(n,\lambda)$:
\[
A=\left\langle z(n,\lambda),z(n,\lambda)+1,z(n,\lambda)+3,\dots, z(n,\lambda)-3,z(n,\lambda)-1\right\rangle.
\]
Define
\[
N(n):=\left\{(f_J\sqcup\langle f_{z(n,\lambda)}\rangle,f_A)\in M(n)\mid
f_J\in\mathbb{F}(n-f_{z(n,\lambda)}-f_{z(n,\lambda-1)})\right\}.
\]
The element $(f_{J_1},f_{J_2})\in M(n)$ belongs to the set
$u\left(\mathbb{F}(n)\right)$
if and only if $(f_{J_1},f_{J_2})\not\in N(n)$.

That is, $u\left(\mathbb{F}(n)\right)=M(n)\setminus N(n)$.
Hence, $\Phi(n;t)$, the generating function of $u\left(\mathbb{F}(n)\right)$,
is the difference of generating functions
of the sets $M(n)$ and $N(n)$.

Since $|A|=a$, it follows that $t^{a+1}\,\Phi\left(n-f_{z(n,\lambda)}-f_{z(n,\lambda-1)};t\right)$
is the generating function of $N(n)$.
\ep

\bp[Proof of Theorem {\rm\ref{th_F}}]
First, let $n$ be an $f$-simple number and $\alpha(n)=(a_1,a_2,\dots, a_\lambda)$.
Let us show that then $\Phi(n;t)=\Delta(a_1,a_2,\dots, a_\lambda;t)$. For $\lambda=1$, this
follows from Lemma \ref{L_Flem02}. This lemma also implies that
$\Phi(f_{z(n,\lambda)-z(n,\lambda-1)+1};t)=\xi(a;t)$. By induction on $\lambda$ we can
rewrite formula \eqref{Leit2} as
\[
\Phi(n;t)=\xi(a_\lambda;t)\,\Delta(a_1,a_2,\dots, a_{\lambda-1};t)
-t^{a_\lambda+1}\,\Delta(a_1,a_2,\dots, a_{\lambda-2};t).
\]
Since the right-hand side of this formula coincides with the right-hand side
of \eqref{eqRec}, we obtain the claim.
In general case,
formula \eqref{PhiD} follows from Lemma \ref{L_Flem05} by induction on $r$.
\ep

Define $F(0)=\chi(0)=1$ and for $n>0$ we set
\[
F(n):=\Phi(n;1)\qquad\text{and}\qquad\chi(n):=\Phi(n;-1).
\]

\bcr\lab{cor1} {\rm (See also \cite{RN}.)}
The values of $\chi(n)$ are $0,\pm 1$ for any $n\in\mathbb{N}$.
\ecr

\bp
It suffices to show that,
for a vector $A=(a_1,\dots, a_m)$ with natural coordinates,
$\Delta(A;-1)=0,\pm 1$. For $m=1,2$, this is clear.

Let $m>2$. Set for brevity $x(m):=\Delta(a_1,\dots, a_m;-1)$.
Then formula \eqref{eqRec} shows that
\[
x(m)=
\bcs
x(m-2)&\text{if $\beta(a_m)=0$},\\
-x(m-1)-x(m-2)&\text{if $\beta(a_m)=1$}.
\ecs
\]
To check that values of $x(m)$ are $0,\pm 1$ only, let us perform induction on $m$.
For $\beta(a_m)=0$, the claim is clear.

For $\beta(a_m)=\beta(a_{m-1})=1$, the claim follows as well, since
\[
x(m)=-x(m-1)-x(m-2)=x(m-2)+x(m-3)-x(m-2)=x(m-3).
\]
If $\beta(a_m)=1$ and $\beta(a_{m-1})=0$, then the claim also follows by induction, because
\[
x(m)=-x(m-1)-x(m-2)=-x(m-3)-x(m-2)=\Delta(a_1,a_2,\dots, a_{m-2},a_{m-1}+1;-1).
\qedhere
\]
\ep

\blr\lab{fi}
We have $n\in[f_i,f_{i+1}-1]$ if and only if $\lambda(n)=i$.
\elr

\bp
The claim obviously follows from the easily verifiable formula
\[
f_{2-\beta(i)}+f_{4-\beta(i)}+{\dots}+f_{i-2}+f_i=f_{i+1}-1.
\qedhere
\]
\ep

This Lemma and Theorem \ref{th_F} imply the next statement, to be used in Sections \ref{chi} and \ref{Fgr}:
\bcr
Let $n\in[f_a,f_{a+1}-1]$, where $0\<a\<i-2$. Then
\beq\lab{PhiTi}
\Phi(n+f_i;t)=
\bcs
\xi\left(\frac{i-a+1}{2};t\right)\cdot\Phi(n;t)&\text{if\; $\beta(a)\neq\beta(i)$},\\[1mm]
\xi\left(\frac{i-a+2}{2};t\right)\cdot\Phi(n;t)-t^{\frac{i-a}{2}+2}\cdot\Phi(n-f_a;t)&\text{if\; $\beta(a)=\beta(i)$}.
\ecs
\eeq
\ecr

The following statement shows that on the interval $[f_i-1,f_{i+1}-1]$ the function $\Phi$ is \quat{symmetric}
in the following sense:

\blr\lab{12}
If $n\in[f_i-1,f_{i+1}-1]$, where $i\>1$, then
\[
\Phi(n;t)=t^i\cdot\Phi\left(f_{i+2}-2-n;\,t^{-1}\right).
\]
In particular, for $t=1$ and for $t=-1$, we, respectively, obtain the formulas:
\beq\lab{fxi}
F(n)=F(f_{i+2}-2-n)\qquad\text{and}\qquad
\chi(n)=(-1)^i\chi(f_{i+2}-2-n).
\eeq
\elr

\bp
Let $\rho_i$ be the reflection of the interval
$[f_i-1,f_{i+1}-1]$ with respect to its center:
\beq\lab{rho}
\rho_i(n)=f_{i+2}-2-n.
\eeq
Since $f_1+f_2+{\dots}+ f_i=f_{i+2}-2$, we see that each Fibonacci partition
$f_I$ of $n\in[f_i-1,f_{i+1}-1]$ with $m$ parts, where $1\<m<i$, corresponds to the Fibonacci partition
$f_{I^\prime}$ of $f_{i+2}-2-n$ with $i-m$ parts, where $I^\prime={\{1,2,\dots, i\}\setminus I}$.
This implies the formula required.
\ep

\section{\bf Solutions of the equation $F(n)=k$ and generating numbers}\lab{Gen}

The main goal of this section is to present an algorithm to
find all solution of the equation $F(n)=k$.
Our main tool here is the notion of generating numbers (Def.\ref{def_GN}).
\bdr
Let $\mathcal{A}$ be the set
of vectors $(a_1,a_2,\dots, a_m)$ with natural coordinates such that
$a_1\>1$ while $a_i\>2$ for all $i\>2$, and let $\mathcal{A}_+\subset\mathcal{A}$
be the set of vectors with $a_1\geqslant 2$.

Any expression of the form $A=A_1\times A_2\times{\dots}\times A_r$, where
$A_1\in\mathcal{A}$ and $A_k\in\mathcal{A}_+$ if $k>1$
is said to be a \emph{polyvector}.
(By definition $A=\emptyset$ if $r=0$.)
The set of polyvectors is denoted by $\mathbb{A}$.
\edr

Obviously, $\mathbb{A}$ coincides with the set
of polyvectors of $2$-partitions (see Def.\ref{defPoly}).
Define
\[
D(\emptyset)=1\qquad\text{and}\qquad D(A)=\Delta(A;1)\quad\text{for $A\in\mathcal{A}$} .
\]
\blr\lab{rac}
Let $A=(a_1,\dots, a_m)\in\mathcal{A}$. Then
\[
\pi(A):=\frac{D(a_2,\dots, a_m)}{D(a_1,a_2,\dots, a_m)}
\]
is an irreducible fraction, $\pi:\mathcal{A}\to\mathbb{Q}_{>0}$ is a one-to-one map,
and $\pi(A)<1$ if and only if $A\in\mathcal{A}_+$.
\elr

It is easy to see that, for $A=(a_1,a_2,\dots, a_m)$, we have
\beq\lab{chain}
\pi(A)=
\cfrac{1}{a_1-
     \cfrac{1}{a_2-
      \cfrac{1}{\ddots-
       \cfrac{1}{a_{m-1}-
        \cfrac{1}{a_m
}}}}}
\eeq
To prove Lemma \ref{rac} it suffices to repeat
(with obvious modifications) the proofs of the corresponding claims on
representations of rational numbers by continued fractions from \cite{HW}, Ch.X.

Let us only note that, in the decomposition
$\frac{l}{k}=\pi(a_1,\dots, a_m)$, the numbers $a_i$ can be calculated
recursively, starting from the values $l_1=l$ and $k_1=k$, by using
the formulas
\beq\lab{frac}
a_i=\left\lceil\frac{k_i}{l_i}\right\rceil,
\qquad l_{i+1}=l_i\cdot a_i-k_i,\qquad k_{i+1}\;=l_i
\qquad\text{if\;\;$l_i\neq 0$}.
\eeq
The calculation terminates when $l_{i+1}=0$.

Define
\[
\G:=\left\{\;g_1\times g_2\times{\dots}\times g_r\;\mid\; g_1\in\mathbb{Q}_{>0}\;\;
\text{and\; $g_i\in\mathbb{Q}_{(0,1)}$\; for $i\>2$;\;\;$r\in\mathbb{N}$}\right\}
\]
and a map $\pi:\mathbb{A}\to\G$ by the formula
\[
\pi(A_1\times{\dots}\times A_r)=\pi(A_1)\times{\dots}\times\pi(A_r).
\]
Lemma \ref{rac} shows that $\pi$ is bijective.
Adding the maps $Z$ and $\alpha$ defined in Section \ref{main}
we obtain a sequence
\beq\lab{chain0}
\mathbb{N}\overset{Z}\longrightarrow\mathbb{P}\overset{\alpha}
\longrightarrow\mathbb{A}\overset{\pi}\longrightarrow\G,
\eeq
where the maps $Z$ and $\pi$ are bijective,
and $\alpha$ is a two-sheeted map.

Theorem \ref{th_F} and Lemma \ref{rac} imply

\bpr\lab{pr33}
Let $\Pi:\mathbb{N}\to\G$ be a map defined by the formula
\[
\Pi(n):=\pi(\alpha(Z(n))).
\]
Define
\[
\G(k):=\left\{\frac{l_1}{k_1}\times\frac{l_2}{k_2}\times{\dots}\times\frac{l_r}{k_r}\in\G\,\mid\,
k_1k_2\dots k_r=k;\; r\in\mathbb{N}\right\}\subset\Gamma.
\]
Then $F(n)=k$ for $n\in\mathbb{N}$ if and only if $\Pi(n)\in\G(k)$.
\epr

The definition of $\alpha$ shows that the sets
\[
\mathbb{P}_a:=\{I\in\mathbb{P}\mid\beta(I^-)=a\},\qquad
\text{where $a=0,1$},
\]
are the sheets of the map $\alpha$ interchanged by the involution
$\sigma:\mathbb{P}\to\mathbb{P}$, defined by the formula
\[
\sigma\langle i_1,i_2,\dots, i_q\rangle:=
\bcs
\langle i_1+1,i_2+1,\dots, i_q+1\rangle &\text{if $\beta(i_1)=1$},\\[0.5mm]
\langle i_1-1,i_2-1,\dots, i_q-1\rangle &\text{if $\beta(i_1)=0$}.
\ecs
\]
In what follows we use the following notation:
\begin{gather*}
\mathbb{N}(k):=\left\{\,n\in\mathbb{N}\,\mid\, F(n)=k\,\right\},\\[1mm]
\mathbb{N}_a:=Z^{-1}(\mathbb{P}_a),\qquad\mathbb{N}_a(k)=\mathbb{N}(k)\cap\mathbb{N}_a,
\qquad\mathbb{P}_a(k)=Z\left(\mathbb{N}_a(k)\right),
\qquad\text{where $a=0,1$}.
\end{gather*}
Then
\[
\mathbb{N}=\mathbb{N}_0\cup\mathbb{N}_1,\qquad \mathbb{N}(k)=\mathbb{N}_0(k)\cup\mathbb{N}_1(k)
\quad\text{and}\quad \mathbb{P}(k)=\mathbb{P}_0(k)\cup \mathbb{P}_1(k).
\]

Since $Z$ is a bijective map, the involution $\sigma:\mathbb{P}\to\mathbb{P}$
induces an involution $\sigma:\mathbb{N}\to\mathbb{N}$
by the formula $\sigma(n):=f_{\sigma(Z(n))}$.
From Theorem \ref{th_F} it follows that $F(n)=F(\sigma(n))$.

Because $\mathbb{N}(k)=\mathbb{N}_1(k)\sqcup\sigma(\mathbb{N}_1(k))$,
to find all solutions of the equation $F(n)=k$
it suffices to compute the set $\mathbb{N}_1(k)$.
For any natural $k$, the sequence \eqref{chain0}
obviously induces the sequence of bijective maps
\beq\lab{chain00}
\mathbb{N}_1(k)\overset{Z}\longrightarrow\mathbb{P}_1(k)\overset{\alpha}
\longrightarrow\mathbb{A}(k)\overset{\pi}\longrightarrow\G(k),
\eeq
where $\mathbb{A}(k):=\pi^{-1}(\G(k))$. To obtain $2$-partition
from $\mathbb{N}_1(k)$ corresponding to $g=g_1\times g_2\times{\dots}\times g_r\in\G(k)$
we can use the next two-step algorithm:
\medskip

\emph{Step 1}: Compute
$A(g)=\pi^{-1}(g_1)\times\pi^{-1}(g_2)\times{\dots}\times\pi^{-1}(g_r)\in\mathbb{A}$ using formulas \eqref{frac}.

\emph{Step 2}: Compute
$\alpha^{-1}(A(g))\in\mathbb{P}_1(k)$.
Then the required Fibonacci partition is $f_{\alpha^{-1}(A(g))}\in\mathbb{N}_1(k)$.
\medskip

For any vector $A=(a_1,a_2,\dots, a_m)\in\mathbb{A}$, define
\begin{align*}
&d(A):=a_1+a_2+{\dots}+ a_m-m,\\
&\upsilon(A):=(a_1+1,a_2+1,\dots, a_m+1),\\
&\varepsilon(A):=\left\langle 2d(a_1)+1,2d(a_1,a_2)+1,\dots, 2d(a_1,a_2,\dots, a_m)+1\right\rangle.
\end{align*}

The following directly verifiable statement computes $\alpha^{-1}(A)\in\mathbb{P}_1$ for any $A\in\mathbb{A}$.

\blr\lab{Eps_L}
Let $A=A_1\times A_2\times{\dots}\times A_r\in\mathbb{A}$.
Then
\beq\lab{eps}
\alpha^{-1}(A)=\varepsilon(A_1)\,\sqcup\,\upsilon^{2d(A_1)+1}\left(\varepsilon(A_2)\right)\,\sqcup\,
{\dots}\,\sqcup\,\upsilon^{2d(A_1)+2d(A_2)+{\dots}+ 2d(A_{r-1})+r-1}\left(\varepsilon(A_r)\right)
\eeq
is the canonical decomposition of $\alpha^{-1}(A)\in\mathbb{P}_1$.

In particular, if $\alpha(Z(n))=A_1\times A_2\times{\dots}\times A_r$, then
\beq\lab{eqz}
z(n,\lambda(n))=2d(A_1)+2d(A_2)+{\dots}+ 2d(A_r)+r.
\eeq
\elr

For example, let us find all solutions of the equation $F(n)=1$.
It is obvious that $\Gamma(1)=\mathbb{N}$. For $a\in\mathbb{N}$, we have
$\pi^{-1}(a)=(1,(a-1)*2)$. From formula \eqref{eps} we obtain
$\alpha^{-1}(\pi^{-1}(a))=\langle 1,3,{\dots},2a-1\rangle\in\mathbb{P}_1(1)$.
Now, formulas
\beq\lab{sum}
f_1+f_3+{\dots}+f_{2a-1}=f_{2a}-1,\qquad
\sigma(f_{2a}-1)=f_2+f_4+{\dots}+f_{2a}=f_{2a+1}-1
\eeq
and Theorem \ref{th_F} show that
\[
\mathbb{N}(1)=\{f_i-1\mid i=1,2,3,{\dots}\}.
\]

The above algorithm certainly is not sufficient to get all solutions of the equation $F(n)=k$
because the set $\G(k)$ is infinite.
Nevertheless, for any $k\>2$, one can canonically define a
\emph{finite set of solutions}, called the $k$-generating numbers,
from which all solutions can be obtained by a regular procedure.
To describe this set of numbers we need some definitions.

Define an action of the (additive) semigroup $\mathbb{Z}_{\>0}$
on the set $\Gamma$ by setting
\[
[a]\left(g_1\times g_2\times{\dots}\times g_r\right)=(g_1+a)\times g_2\times{\dots}\times g_r,
\]
where $a\in\mathbb{Z}_{\>0}$.
Define an involution $\tau:\Gamma\to\Gamma$ by the formulas
\begin{align*}
\tau(l)&=\;l,\\
\tau(g_1\times g_2\times\dots\times g_r)&=\;1\times g_1\times g_2\times\dots\times g_r,\\
\tau(l\times g_1\times g_2\times\dots\times g_r)&=\;(l-1+g_1)\times g_2\times\dots\times g_r,\\
\tau((l+g_1)\times g_2\times\dots\times g_r)&=\;(l+1)\times g_1\times g_2\times\dots\times g_r,
\end{align*}
where $l\in\mathbb{N}$, $g_1,g_2,{\dots},g_r\in\mathbb{Q}_{(0,1)}$ and $r\>1$.
Let $T(\tau)$ be the group with two elements generated by $\tau$.

It is trivial to check that $[a]\cdot\tau=\tau\cdot[a]$.
Therefore, an action of the semigroup $H_0:=\mathbb{Z}_{\>0}\times T(\tau)$ on $\Gamma$ is canonically defined;
obviously, $H_0(\Gamma(k))\subset\Gamma(k)$.
The next claim easily follows from the definitions.

\blr\lab{fund}
The canonical action of the semigroup $H_0$ on the set
$\G^*:=\Gamma\setminus\Gamma(1)$ is free and
\[
\G_+:=\left\{\text{$g_1\times g_2\times{\dots}\times g_r\in\G\mid g_1,g_2,{\dots},g_r\in\mathbb{Q}_{(0,1)}$}\right\}\subset\G^*
\]
is a fundamental domain of this action.
That is, $H_0(\Gamma_+)=\G^*$
and if $g\in\Gamma_+$, then $h(g)\not\in\Gamma_+$ for any nonunit element $h\in H_0$.
\elr

Using the sequence of the bijective maps
\[
\mathbb{N}_1\overset{Z}\longrightarrow
\mathbb{P}_1\overset{\alpha}\longrightarrow
\mathbb{A}\overset{\pi}\longrightarrow\Gamma,
\]
we transfer the action of $[a]$ and $\tau$ on $\mathbb{N}_1$.
Further, using the decomposition $\mathbb{N}=\mathbb{N}_0\sqcup\mathbb{N}_1$
and the involution $\sigma:\mathbb{N}\to\mathbb{N}$ we extend the actions of $[a]$ and $\tau$
on $\mathbb{N}$ by the formulas
\begin{align*}
[a](n)=&\;\sigma\cdot[a]\cdot\sigma(n),\\
\tau(n)=&\;\sigma\cdot\tau\cdot\sigma(n),
\end{align*}
where $n\in\mathbb{N}_0$. Thus, we see that an action of the semigroup
\[
H:=\mathbb{Z}_{\>0}\times T(\tau)\times T(\sigma)
\]
on the set $\mathbb{N}$ is canonically defined. Here $T(\tau)$ and $T(\sigma)$ are the groups,
each with two elements, generated by the involutions $\tau$ and $\sigma$, respectively.

\bdr\lab{def_GN}
Any number $n\in\mathbb{N}_1$, such that $z(n,1)\>3$, is said to be a \emph{generating number}.
If, in addition, $F(n)=k$, then $n$ is said to be a \emph{$k$-generating number}.

By $\mathbb{G}$ and $\mathbb{G}(k)$ we denote the set of all generating
and the set of $k$-generating numbers, respectively.
\edr

For example, $f_{2k-1}$ is a $k$-generating number for $k>1$ since $F(f_{2k-1})=k$;
so $\mathbb{G}(k)\neq\emptyset$ for any $k>1$.
The first 20 generating numbers $n$ and the corresponding values of $F(n)$ are as follows (see \cite{SL},\,A003623):
\smallskip
\begin{center}
\begin{tabular}{|c||c|c|c|c|c|c|c|c|c|c|c|c|c|c|c|c|c|c|c|c|c|c|c|c|c|}\hline
$n$ & 3 & 8 & 11 & 16 & 21 & 24 & 29 & 32 & 37 & 42 & 45 & 50 & 55 & 58 & 63 & 66 & 71 & 76 & 79 & 84\\
\hline
$F(n)$ & 2 & 3 & 3 & 4 & 4 & 5 & 5 & 4 & 6 & 6 & 6 & 6 & 5 & 7 & 8 & 7 & 8 & 7 & 8 & 7\\
\hline
\end{tabular}
\end{center}
\smallskip

The definitions imply
\blr\lab{fund0}
A number $n\in\mathbb{N}_1$ is a generating number whenever $\Pi(n)\in\Gamma_+$, and
$n\in\mathbb{N}_1$ is a $k$-generating number whenever $\Pi(n)\in\Gamma_+(k):=\G(k)\bigcap\G_+$.
In particular, the map $\Pi:\mathbb{G}(k)\to\G_+(k)$ is bijective.
\elr

Since $\G_+(k)$ is a finite set,
$\mathbb{G}(k)$ is a finite set as well. Define
\[
\Psi(k):=\left|\G_+(k)\right|=\left|\mathbb{G}(k)\right|.
\]
For some details on the function $\Psi(k)$, see Section \ref{Add}.

The set $\G_+$ has a structure of a monoid freely generated by the set $\mathbb{Q}_{(0,1)}$
with the multiplication $\times$. It is obvious that
$\G_+(k_1)\times\G_+(k_2)\subset\G_+(k_1k_2)$. Therefore, Lemma \ref{fund0} implies that
the set of generating numbers $\mathbb{G}=\coprod_{k\>2}\mathbb{G}(k)\subset\mathbb{N}$ also has a canonical structure
of a monoid isomorphic to $\G_+$. Denote the multiplication in $\mathbb{G}$
by $\times$ as well. Then
$\mathbb{G}(k_1)\times\mathbb{G}(k_2)\subset\mathbb{G}(k_1k_2)$, and
the map $F:\mathbb{G}\to\mathbb{N}$ is a homomorphism of monoids.

In other words, we can $\times$-multiply the generating solutions
$n_1$ and $n_2$ of the equations $F(n)=k_1$ and $F(n)=k_2$ and
obtain, for $n_1\neq n_2$, two different generating solutions
$n_1\times n_2$ and $n_2\times n_1$ of the equation $F(n)=k_1k_2$.
For example, $F(3)=2$ and $F(8)=3$. Then $3\times 8=37$ and
$8\times 3=42$ are generating solutions of the equation $F(n)=6$.

Lemma \ref{fund} and Lemma \ref{fund0} directly imply

\btr\lab{Gplus}
The semigroup $H\simeq\mathbb{Z}_{\>0}\times\mathbb{Z}_2\times\mathbb{Z}_2$
canonically acts on the set $\mathbb{N}$ so that,
for any $h\in H$, we have $F(h(n))=F(n)$.

The restriction of the $H$-action on the set $\mathbb{N}\setminus\mathbb{N}(1)$ is free,
and the set $\mathbb{G}\subset\mathbb{N}\setminus\mathbb{N}(1)$
of generating numbers is a fundamental domain of this restricted action.
In particular, for any $k\>2$, we have $\mathbb{N}(k)=H(\mathbb{G}(k))$.
\etr

It is convenient to describe the $H$-action on $n\in\mathbb{N}\setminus\mathbb{N}(1)$ explicitly.
Define
\[
\theta(n):=
\bcs
1&\text{if\; $z(n,1)\>3$},\\
\min\{\,m\mid z(n,m)-z(n,m-1)\>3\,\}&\text{if\; $z(n,1)<3$\; and\; $1<m\<\lambda(n)$},
\ecs
\]
\[
Z^\prime(n):=\langle z(n,\theta(n)),z(n,\theta(n)+1),{\dots},z(n,\lambda(n))\rangle.
\]
For any partition $I=\langle i_1,i_2,\dots,i_q\rangle$
and any integer $r$, we set $I+r:=\langle i_1+r,i_2+r,\dots,i_q+r\rangle$.
Thanks to formulas \eqref{sum} we have (see formula \eqref{len} for notation)
\beq\lab{nsig}
n=f_{2\theta(n)-\beta(z(n,1))-1}-1+\|f_{Z^\prime(n)}\|
\quad\text{and}\quad
\sigma(n)=
\bcs
f_{2\theta(n)-1}+\|f_{Z^\prime(n)+1}\|&\text{if $\beta(z(n,1))=1$},\\
f_{2\theta(n)-2}+\|f_{Z^\prime(n)-1}\|&\text{if $\beta(z(n,1))=0$}.
\ecs
\eeq
If $n\in\mathbb{N}_1\setminus\mathbb{N}_1(1)$,
then from the definitions we obtain formulas
\beq\lab{Sma}
[a](n)=f_{2\theta(n)-2+2a}-1+\|f_{Z^\prime(n)+2a}\|,
\eeq
\beq\lab{Smtau}
\tau(n)=
\bcs
f_{2\theta(n)}-1+\|f_{Z^\prime(n)+1}\|&\text{if\; $\beta(z(n,\theta(n)))=1$},\\
f_{2\theta(n)-4}-1+\|f_{Z^\prime(n)-1}\|&\text{if\; $\beta(z(n,\theta(n)))=0$},
\ecs
\eeq
Since $[a]\cdot\sigma=\sigma\cdot[a]$ and $\tau\cdot\sigma=\sigma\cdot\tau$,
formulas \eqref{nsig}-\eqref{Smtau} completely describe the $H$-action
on $\mathbb{N}\setminus\mathbb{N}(1)$.
\smallskip

\noindent{\bf Examples.}\; Let us apply Theorem \ref{Gplus} to find all solutions of the
equation\footnote{The solutions of the equation $F(n)=k$ for $k=1,2,3$
were earlier obtained in \cite{KL}}
$F(n)=k$, where $k=2,3,4$. In these examples, to compute the set
$\mathbb{G}(k)$ by means of the set $\Gamma_+(k)$ we use formula \eqref{eps}.
\smallskip
\newline
Since $\G_+(2)=\left\{\frac{1}{2}\right\}$, we have $\mathbb{G}(2)=\{f_3\}$.
Then formulas \eqref{Sma} and \eqref{Smtau} show that any
solution $n\in\mathbb{N}_1$ of the equation $F(n)=2$
can be uniquely represented in one of the forms:
\[
\text{either}\qquad[a](f_3)=f_{2a}-1+f_{2a+3},\qquad\text{or}\qquad
[a](\tau(f_3))=[a](f_2-1+f_4)=f_{2(a+1)}-1+f_{2a+4},
\]
where $a\>0$.
Applying involution $\sigma$ to these numbers
we see that the set of solutions of the equation $F(n)=2$
coincides with the set of numbers
\[
f_i-1+f_{i+3},\quad f_{i+2}-1+f_{i+4},\qquad\text{where\; $i\in\mathbb{Z}_{\>0}$}.
\]
Let $k=3$ or $k=4$. Because $\G_+(3)=\left\{\frac{1}{3},\frac{2}{3}\right\}$
and $\G_+(4)=\left\{\frac{1}{4},\frac{3}{4},\frac{1}{2}\times\frac{1}{2}\right\}$,
we, respectively, have $\mathbb{G}(3)=\{f_5,f_3+f_5\}$ and
$\mathbb{G}(4)=\left\{f_7,f_3+f_5+f_7,f_3+f_6\right\}$.
Repeating the similar procedure as for $k=2$, we see that the set of solutions
of the equation $F(n)=3$ coincides with the set of numbers
\[
f_i-1+f_{i+5},\quad f_{i+2}-1+f_{i+6},\quad f_i-1+f_{i+3}+f_{i+5},\quad
f_{i+2}-1+f_{i+4}+f_{i+6},
\]
and the set of solutions of the equation $F(n)=4$ coincides with the set of numbers
\begin{gather*}
f_i-1+f_{i+7},\quad f_{i+2}-1+f_{i+8},\quad
f_i-1+f_{i+3}+f_{i+5}+f_{i+7},\quad f_{i+2}-1+f_{i+4}+f_{i+6}+f_{i+8},\\[1mm]
f_i-1+f_{i+3}+f_{i+6},\quad f_{i+2}-1+f_{i+4}+f_{i+7},
\end{gather*}
where $i\in\mathbb{Z}_{\>0}$.
\smallskip

The set of generating numbers has the following compact description:

\btr\lab{G}
The natural number $n$ is a generating number if and only if $n=2\lb l\phi\rb+l$,
where $l\in\mathbb{N}$.
\etr

\bp
Using the known formula
$\lb l\phi\rb =f_{z(l,1)+1}+{\dots}+ f_{z(l,q)+1}-\beta(z(l,1))$ (see \cite{GKP}, p. 339),
we obtain
\[
2\lb l\phi\rb +l=f_{z(l,1)+3}+{\dots}+ f_{z(l,q)+3}-2\beta(z(l,1)).
\]
Since $f_{z(l,1)+3}-2=f_3+f_5+{\dots}+ f_{z(l,1)+2}$, we have
\beq\lab{eqess}
2\lb l\phi\rb+l=
\bcs
f_{z(l,1)+3}+{\dots}+f_{z(l,q)+3} &\text{if}\quad\beta(z(l,1))=0,\\[1mm]
f_3+f_5+{\dots}+f_{z(l,1)+2}+f_{z(l,2)+3}+{\dots}+f_{z(l,q)+3} &\text{if}\quad\beta(z(l,1))=1.
\ecs
\eeq
Hence $2\lb l\phi\rb+l\in\mathcal{\mathbb{G}}$.

Conversely, let $n\in\mathbb{G}$. Set
\[
u(n)=
\bcs
\left\langle z(n,1)-3,z(n,2)-3,\dots,z(n,\lambda(n))-3\right\rangle &\text{if $z(n,1)>3$},\\[1mm]
\left\langle 2a-1,z(n,a+1)-3,\dots,z(n,\lambda(n))-3\right\rangle &\text{if $Z(n)=\langle 3,5,\dots, 2a+1,z(n,a+1),\dots, z(n,\lambda(n))\rangle$},
\ecs
\]
where $z(n,a+1)>2a+3$. Then, for $l=f_{u(n)}$, we have $n=2\lb l\phi\rb+l$
by formula \eqref{eqess}.
\ep

\section{\bf Estimates of $D(A)$}\lab{EstD}

Here we establish an auxiliary statement to be used in the next two sections.
The proof of it is based on the following formal identities:
\begin{gather}
D(a_1,\dots, a_s,k*2)=D(a_1,\dots, a_s)+kD(a_1,\dots, a_{s-1},a_s-1),
\lab{T1}\\[1mm]
D(a_1,\dots, a_m)=D(a_1,\dots, a_{m-1},2)+ (a_m-2)D(a_1,\dots, a_{m-1}).\lab{T2}
\end{gather}

\blr\lab{L_lowD}
Let $A=(a_1,\dots, a_m)\in\mathbb{A}_+$. Then
\beq\lab{dD}
d(A)+1\<D(A)\<f_{d(A)+1}.
\eeq
\noindent
The left inequality in \eqref{dD} becomes an equality if and only if $m=1$ or $A=(2,2,{\dots},2)$.
\newline
\noindent
The right inequality in \eqref{dD} becomes an  equality if and only if $a_1\<3$,
$i_m\<3$ and $a_2={\dots}=a_{m-1}=3$.
\elr

\bp
First, consider the left inequality. For $m=1$ or $A=(2,2,\dots,
2)$, the claim is clear. For $m>1$, let us perform induction on $m$. To begin with,
assume that $A=(a_1,\dots, a_s,k*2)$, where $a_s>2$.

Since $a_1\>2$ and $a_s-1\>2$, then formula
\eqref{T1} and the inductive hypothesis imply that
\[
D(A)>(k+1)(a_1+{\dots}+ a_s-s)+1>a_1+{\dots}+ a_s+k-s+1=d(A)+1.
\]
Let $a_m>2$. Formula \eqref{T2},
the inductive hypothesis, and the already considered case $a_m=2$ show that
\[
D(A)>(a_m-1)(a_1+{\dots}+ a_{m-1}-m+2)+1>a_1+{\dots}+ a_m-m+1=d(A)+1.
\]
The left inequality is proved.

To prove the right inequality we use an obvious claim:
\beq\lab{triv} f_{a+1}+lf_a\<f_{a+l+1},\quad\text
{where equality holds if and only if\; $l=0,1$.}
\eeq

If $A=(2,2,\dots, 2)$, then $D(A)=d(A)+1\<f_{d(A)+1}$,
where the equality holds whenever $A=(2)$ or $A=(2,2)$.
Otherwise we use the induction on $m$.

Let $A=(a_1,\dots, a_s,k*2)$, where $a_s>2$.
Then formula \eqref{T1} and the inductive hypothesis show
that the inequality of lemma is a corollary of the inequality
\eqref{triv} when $a=a_1+a_2+{\dots}+ a_s-s$ and $l=k$; this inequality
turns into an equality only for $k=1$. In this case, the inductive
hypothesis implies also that $a_1\<3$ and $a_2={\dots}=a_{m-1}=3$.

Let $a_m>2$. Then formula \eqref{T2}, the inductive hypothesis,
and the already considered case $a_m=2$ show that the right inequality
in \eqref{dD} is a corollary of inequality \eqref{triv} for
$a=a_1+{\dots}+ a_{m-1}-m+2$ and $s=a_m-2$.
Moreover, the inductive hypothesis also
implies that the right inequality in \eqref{dD} turns into an equality only when
$a_1\<3$ and $a_2={\dots}=a_m=3$.
\ep

\section{\bf Stabilization of Fibonacci partitions}\lab{Stab}

In this section we are using the results of Sections \ref{Gen} and \ref{EstD} to
show that the function $F(n)$ stabilizes in the following sense:
\btr\lab{2Psi}
The quantity of numbers $n\in\mathbb{N}$ such that
$f_i-1\<n<f_{i+1}-1$ and $F(n)=k$ does not depend on $i$ for $i\>2k$.
It is equal to  $2\Psi(k)$ if $k>1$, and to $1$ if $k=1$.
\etr

\blr\lab{Ups}
Let $n\in\mathbb{G}$ and
$\alpha(Z(n))=A_1\times A_2\times{\dots}\times A_r$. Then
$n\<f_{2d(A_1)+f_{2d(A_2)}+{\dots}+ 2d(A_r)+r+1}-2$.
\elr

\bp
Formula \eqref{eqz} shows that $z(n,\lambda(n))=2d(A_1)+2d(A_2)+{\dots}+ 2d(A_r)+r$.
Since $n\in\mathbb{G}$ we have
\[
n=f_{z(n,1)}+f_{z(n,2)}+{\dots}+f_{z(n,\lambda(n))}\<f_3+f_5+{\dots}+f_{z(n,\lambda(n))}=f_{2d(A_1)+2d(A_2)+{\dots}+ 2d(A_r)+r+1}-2.
\qedhere
\]
\ep

\blr\lab{press}
The maximal $k$-generating number is equal to $f_{2k}-2$.
\elr

\bp
First, note that $f_{2k}-2$ is a $k$-generating number since
\[
\Pi(f_{2k}-2)=\Pi\left(f_3+f_5+{\dots}+f_{2k-1}\right)=\frac{k-1}{k}\in\G_+.
\]

Let $n$ be the maximal $k$-generating number and $\alpha\left(Z(n)\right)=A_1\times{\dots}\times A_r$.
From Theorem \ref{th_F} and Lemma \ref{L_lowD} we obtain
\[
k=D(A_1)\cdot{\dots}\cdot D(A_r)\>(d(A_1)+1)\cdot{\dots}\cdot(d(A_r)+1)\>d(A_1)+{\dots}+d(A_r)+r.
\]
This inequality and Lemma \ref{Ups} show that
\[
n\<f_{2d(A_1)+{\dots}+ 2d(A_r)+r+1}-2\<f_{2k-r+1}-2\<f_{2k}-2.
\]
Since $f_{2k}-2$ is a $k$-generating number, it follows
that $n=f_{2k}-2$ in view of the choice of $n$.
\ep

\bp[Proof of Theorem {\rm\ref{2Psi}}]
For $k=1$, the statement follows from the above description of $\mathbb{N}(1)$.
Let $k\>2$.

In the proof we use the action
of the semigroup $H$ on $\mathbb{N}\setminus\mathbb{N}(1)$,
defined in Section \ref{Gen}. Remind that $H$ is generated
by the elements $[a]$, where $a\in\mathbb{Z}_{\>0}$,
and by the involutions $\tau$ and $\sigma$;
see formulas \eqref{nsig}--\eqref{Smtau}
for their action on $\mathbb{N}\setminus\mathbb{N}(1)$.
For brevity, set
\[
V(m):=[f_{m}-1,f_{m+1}-1].
\]
Lemma \ref{fi}, Lemma \ref{press}, and formulas \eqref{Sma}, \eqref{Smtau}
imply, that for any $n\in\mathbb{G}(k)$ and any $m\>2k-1$,
there exists a unique $a=a(n,m)\in\mathbb{Z}_{\>0}$ such that
either number $[a](n)$, or number $[a](\tau(n))$
belongs to the set $V(m)$.
Let $B_k(m)$ be the set of such numbers obtained for all $n\in\mathbb{G}(k)$.
From Theorem \ref{Gplus} it follows that
\[
B_k(m)=\mathbb{N}_1(k)\cap V(m).
\]
Then $\mathbb{N}_0(k)\cap V(m+1)=\sigma(B_k(m))$.
Thus, for any $m\>2k$, we have
\[
\mathbb{N}(k)\cap V(m)=\sigma(B_k(m-1))\sqcup B_k(m),\qquad\text{where}\qquad
\sigma(B_k(m-1))\cap B_k(m)=\emptyset.
\]
Since $|B_k(m)|=|\mathbb{G}(k)|$ for $m\>2k-1$, then
$|\mathbb{N}(k)\cap V(m)|=2\Psi(k)$ for any $m\>2k$, as claimed.
\ep

\section{\bf The upper bounds for $F(n)$}\lab{Up}

In this section, for any $i\in\mathbb{N}$, we obtain the least upper bound for $F(n)$
when $n\in[f_i,f_{i+1}-1]$.

\blr\lab{ineq_lm}
Let $i_1,i_2,{\dots} ,i_r$ be the natural numbers $\>2$. Then
\beq\lab{ineq}
f_{i_1}\cdot f_{i_2}\cdot{\dots}\cdot f_{i_r}\<
\bcs
f_{i_1+i_2+{\dots}+i_r-\frac{r-1}{2}}&\text{\rm if\; $\beta(r)=1$},\\[1mm]
2\,f_{i_1+i_2+{\dots}+i_r-\frac{r}{2}-1}&\text{\rm if\; $\beta(r)=0$}.
\ecs \eeq
These inequalities turn into equalities only if either $r=1$,
or one of the following cases takes place:

$r=2$ and either $i_1=2$, or $i_2=2$;

$r=3$ and $i_1=i_2=i_3=2$;

$r=4$ and $i_1=i_2=i_3=i_4=2$.
\elr

\bp
The well known formula (see \cite{GKP}, Sec.6.6)
\begin{equation}\lab{St3}
f_{a+b}=f_a f_b+f_{a-1}f_{b-1},\text{~~where $a,b\in\mathbb{N}$,}
\end{equation}
easily implies the inequalities
\begin{equation}\lab{St03}
f_{a+b-1}\<f_af_b\<\frac12f_{a+b+1}\,.
\end{equation}
The left inequality in \eqref{St03} turns into an equality only if either $a=1$, or
$b=1$, and the right inequality in \eqref{St03} turns into an equality only when
$a=b=2$.

To prove Lemma \ref{ineq_lm} let us perform induction on $r$.
For $r=1$, the claim is trivial.

Applying formula \eqref{St3} to $a=i_1$, $b=i_2$ and to $a=2$,
$b=i_1+i_2-2$, we see that
\begin{equation}\lab{*}
f_{i_1}f_{i_2}=2f_{i_1+i_2-2}-(f_{i_1-1}f_{i_2-1}-f_{i_1+i_2-3})\<2f_{i_1+i_2-2}.
\end{equation} The inequality \eqref{*} follows from the left
inequality in \eqref{St03}, where $a=i_1-1$ and $b=i_2-1$; and
\eqref{*} turns into an equality only if either $i_1=2$, or $i_2=2$.
This proves Lemma for $r=2$.

Let $r>2$. For $\beta(r)=0$, the inductive hypothesis gives the
inequality required, since
\[
f_{i_1}\cdot{\dots}\cdot f_{i_{r-1}}\cdot f_{i_r}\<f_{i_1+{\dots}+i_{r-1}-\frac{r-2}{2}}\cdot f_{i_r}
\<2f_{i_1+{\dots}+i_r-\frac{r}{2}-1}.
\]
For $\beta(r)=1$, the inequality required implied by the inductive
hypothesis  and the right inequality in \eqref{St03}:
\[
f_{i_1}\cdot{\dots}\cdot f_{i_{r-1}}\cdot f_{i_r}\<2f_{i_1+{\dots}+i_{r-1}-\frac{r+1}{2}}\cdot f_{i_r}\<
f_{i_1+{\dots}+i_r-\frac{r-1}{2}}.
\]
The equalities in \eqref{ineq} are also covered by the last two
inequalities.
\ep

\btr\lab{Up_Est_r}
Let $n\in[f_i,f_{i+1}-1]$ and $i\>3$. Then
\beq
F(n)\<
\bcs\lab{thest}
f_{\frac{i+1}{2}}&\text{\rm if\; $\beta(i)=1$},\\[1mm]
2f_{\frac{i}{2}-1}&\text{\rm if\; $\beta(i)=0$ and $i\neq 2$},\\[1mm]
1&\text{\rm if\; $i=2$}. \ecs \eeq Let $M(i)\subset[f_i,f_{i+1}-1]$ be
the set of numbers $n$ for which inequalities \eqref{thest} become
equalities.
For $i\>5$, define the set of numbers
\[
M_i:=\\
\bcs
\left\{\,f^2_{\frac{i+1}{2}}-1,\;f^2_{\frac{i-1}{2}}+f_i-1\,\right\}
&\text{\rm if\; $\beta(i)=1$},\\[2mm]
\Big\{\,
f_{\frac{i}{2}}f_{\frac{i}{2}+1}+2\cdot(-1)^{\frac{i}{2}+1}-1,\;
f^2_{\frac{i}{2}-1}+f_i-1,\;\\
\qquad f^2_{\frac{i}{2}-2}+2f_{i-1}-1,\;
f_{\frac{i}{2}-1}f_{\frac{i}{2}}+f_i+2\cdot(-1)^{\frac{i}{2}}-1
\,\Big\}
&\text{\rm if\; $\beta(i)=0$}.
\ecs
\]
Then $M(i)=M_i$ if $i\not\in\{1,2,3,4,6,9,12\}$, and
\begin{gather*}
M(1)=\{1\},\qquad M(2)=\{2\},\qquad M(3)=\{3\},\qquad M(4)=\{5,6\},\\
M(6)=\{16\},\qquad M(9)=M_9\cup\{71\},\qquad M(12)=M_{12}\cup\{304\}.
\end{gather*}
\etr

\bp
We say that vector $A=(a_1,a_2,{\dots},a_q)\in\mathcal{A}_+$ is \textit{special}
if $a_1\<3$ and $a_2={\dots}=a_q=3$. Since,
for any special vector, we have
\[
d(A)=
\bcs
2q-1&\text{if\; $A=\left(2,(q-1)*3\right)$},\\[1mm]
2q&\text{if\; $A=(q*3)$}.
\ecs
\]
it follows that for any vector $A\in\mathcal{A}_+$, there exists a unique
special vector $\widehat{A}$ such that $d(A)=d(\widehat{A})$.
Let
\[
\alpha(Z(n))=A_1\times A_2\times{\dots}\times A_r,\qquad
d(n)=d(A_1)+d(A_2)+{\dots}+d(A_r),\qquad
r(n)=r.
\]
By formula \eqref{eqz}
\beq\lab{nri}
n\in[f_i,f_{i+1}-1]\qquad\text{if and only if\qquad $2d(n)+r(n)=i$}.
\eeq

Let
$\widehat{n}:=Z^{-1}\alpha^{-1}(\widehat{A}_1\times\widehat{A}_2\times{\dots}\times\widehat{A}_r)$.
Since $d(\widehat{n})=d(n)$  and $r(\widehat{n})=r(n)$, it follows
that $\widehat{n}\in[f_i,f_{i+1}-1]$. Moreover, by the right
inequality in \eqref{dD} and Theorem \ref{th_F} we have
\begin{equation}\lab{product}
F(n)\<F(\widehat{n})=f_{d(A_1)+1}\cdot f_{d(A_2)+1}\cdot{\dots}\cdot f_{d(A_r)+1}.
\end{equation}
Thus, we can assume that $n=\widehat{n}$, i.e., all vectors $A_k$
are special.

Let $A(i)$ be a special vector such that (use formula \eqref{nri})
\beq\lab{Ai}
d(A(i))=
\bcs
d(A_1)+d(A_2)+{\dots}+d(A_r)+\frac{r-1}{2}=\frac{i-1}{2}&\text{\rm if\; $\beta(i)=1$},\\[1mm]
d(A_1)+d(A_2)+{\dots}+d(A_r)+\frac{r}{2}-2=\frac{i-4}{2}&\text{\rm if\; $\beta(i)=0$}.
\ecs
\eeq

Define the number $n(i)\in[f_i,f_{i+1}-1]$ by the formula
\[
n(i):=
\bcs
Z^{-1}\alpha^{-1}(\,A(i)\,)&\text{\rm if\; $\beta(i)=1$},\\[1mm]
Z^{-1}\alpha^{-1}(\,(2)\times A(i)\,)&\text{\rm if\; $\beta(i)=0$}.
\ecs
\]
Let us show that $F(n)\<F(n(i))$.
Indeed, the assumption on $n$ and formulas  \eqref{product} and \eqref{Ai}
imply that this inequality is equivalent to the inequality
\[
f_{d(A_1)+1}\cdot f_{d(A_2)+1}\cdot{\dots}\cdot f_{d(A_r)+1}\<f_{d(A(i))+1}=
\bcs
f_{d(A_1)+d(A_2)+{\dots}+d(A_r)+\frac{r-1}{2}+1}=f_{\frac{i+1}{2}}&\text{\rm if\; $\beta(i)=1$},\\[1mm]
2\,f_{d(A_1)+d(A_2)+{\dots}+d(A_r)+\frac{r}{2}-1}=2f_{\frac{i}{2}-1}&\text{\rm if\; $\beta(i)=0$}.
\ecs
\]
The substitution $d(A_k)=i_k-1$ turns it into the inequality \eqref{ineq}.
Thus, \eqref{thest} is proved.

The same reasoning and Lemma \ref{L_lowD} show that for any $n\in M(i)$,
we have
\[
\alpha(Z(n))=
\bcs
A&\text{if $\beta(i)=1$},\\
(2)\times A\;\;\text{or}\;\;A\times(2)&\text{if $\beta(i)=0$},
\ecs
\]
where $A$ is of one of the following forms
\[(2,3,{\dots},3),\quad (3,{\dots},3,2),\quad
(3,3,{\dots},3),\quad (2,3,{\dots},3,2),\quad (2)\times(2)\times(2).
\]

Namely, formula \eqref{nri} implies that if $\beta(i)=1$, then
\[
A=
\bcs
(2)\times(2)\times(2)&\text{if $i=9$},\\
(2)\times(2)\times(2)\times(2)&\text{if $i=12$}
\ecs
\]
If $i\neq 9,12$ and if $\beta(i)=1$, then
\[
A=
\bcs
\left(2,\frac{i-3}{4}*3\right)\;\;\text{or}\;\;\left(\frac{i-3}{4}*3,2\right)&\text{if $i\equiv 3\mod 4$},\\[1mm]
\left(\frac{i-1}{4}*3\right)\;\;\text{or}\;\;\left(2,\frac{i-5}{4}*3,2\right)&\text{if $i\equiv 1\mod 4$},
\ecs
\]
and if $\beta(i)=0$, then
\[
A=
\bcs
\left(2,\frac{i-6}{4}*3\right)\;\;\text{or}\;\;\left(\frac{i-6}{4}*3,2\right)&\text{if $i\equiv 2\mod 4$},\\[1mm]
\left(\frac{i-4}{4}*3\right)\;\;\text{or}\;\;\left(2,\frac{i-8}{4}*3,2\right)&\text{if $i\equiv 0\mod 4$},
\ecs
\]

Since the polyvectors $\alpha(Z(n))$ for all $n\in M(i)$
are now known, we can find the corresponding $n$
by a routine calculation, thanks to formula \eqref{eps} and easily
verifiable identities
\begin{align*}
&f_3+f_7+{\dots}+f_{4q-1}=f^2_{2q}-1,
&&f_5+f_9+{\dots}+f_{4q+1}= f^2_{2q+1}-1,\\[1mm]
&f_2+f_6+{\dots}+f_{4q-2}=f_{2q-1}\cdot f_{2q},
&&f_4+f_8+{\dots}+f_{4q}= f_{2q}\cdot f_{2q+1}-1.
\end{align*}
As a result, we obtain the claimed formulas for $M_i$.

When $i=9$ or $i=12$, we should add $f_3+f_6+f_9=71$ and
$f_3+f_6+f_9+f_{12}=304$, respectively, to the numbers constituting the sets $M_i$.
This completes the proof.
\ep

\bcr
If $i\not\in\{1,2,3,4,6,9,12\}$, then
\[
\left|M(i)\right|=
\bcs
2&\text{\rm if $\beta(i)=1$},\\[1mm]
4&\text{\rm if $\beta(i)=0$}.
\ecs
\]
\ecr
\bcr\lab{Up_Est0}
For $n\in\mathbb{N}$, we have
$F(n)\<\sqrt{n+1}$, where the equality holds whenever $n=f_a^2-1$.
\ecr

\bp Since the function $\sqrt{n+1}$ grows together with $n$, it
suffices to check that $F(n)\<\sqrt{n+1}$ only for the minimal
numbers $n\in M(i)$ and arbitrary $i\>1$. A direct calculation shows
that for $i\>5$, we have
\beq\lab{minM}
\min\,M(i)=
\bcs
f^2_{\frac{i+1}{2}}-1&\text{if $\beta(i)=1$},\\[1.5mm]
f_{\frac{i}{2}}f_{\frac{i}{2}+1}+2\cdot(-1)^{\frac{i}{2}+1}-1&\text{if $\beta(i)=0$}.
\ecs
\eeq
Thus, for $\beta(i)=1$, the claim directly follows from Theorem
\ref{Up_Est_r}.

Let $\beta(i)=0$ and $i=2m$.  For $m=1,2,3,4$, the claim is trivial.
Thanks to the formula \eqref{minM} and Theorem \ref{Up_Est_r} to finish the
proof it suffices to show that
\[
4f^2_{m-1}<f_mf_{m+1}+2\cdot(-1)^{m+1}\text{~~ for $m\>5$}.
\]
Using the \textit{Cassini formula}
$f^2_{m-1}=f_mf_{m-2}+(-1)^{m+1}$ (see \cite{GKP}, Sec.6.6), we obtain
\[
f_mf_{m+1}-4f^2_{m-1}+2\cdot(-1)^{m+1}=
f_m\left(f_{m+1}-4f_{m-2}\right)+2(-1)^m=f_mf_{m-5}+2(-1)^m>0.
\qedhere
\]
\ep

\section{\bf How often does $\chi(n)$ vanish}\lab{chi}

In this section we use Theorem \ref{th_F} to prove the next statement --
the main result of paper \cite{ARP}:
\btr\lab{th41}
Let $E(a):=\left|\{n\mid 0\<n\<a\;\;\text{and}\;\;\chi(n)=0\}\right|$. Then
$\mathop{\lim}\limits_{a\to\infty}\frac{E(a)}{a}=1$.
\etr

First we establish several auxiliary formulas.

\blr\lab{14}
Let $n\in[f_{i-3}+f_i-1,f_{i-2}+f_i-1]$, where $i\>3$.
Then\footnote{This claim appeared also in \cite{RN} in the proof of Theorem 9.} $\chi(n)=0$.
\elr

\bp
Lemma follows from formula \eqref{PhiTi}, where $t=-1$ and $a=i-3$.
\ep

\blr\lab{15}
Let $n\in[0,f_i-1]$ and $s\>i\>1$. Then $\chi(n)=\chi(n+f_s+f_{s+2})$.
\elr

\bp
It suffices to assume that $n$ is an $f$-simple number.
Let $Z(n)=\langle i_1,i_2,\dots,i_m\rangle$ and
$\alpha(n)=(a_1,a_2,\dots, a_m)$.

By Lemma \ref{L_Flem01} we have $i_m\<i-1$. If $i_m=i-1$ and $s=i$, then

\[
n+f_i+f_{i+2}=f_{i_1}+\dots+f_{i_{m-1}}+f_{i+3}.
\]
Therefore, $\alpha(n+f_i+f_{i+2})=(a_1,\dots, a_{m-1},a_m+2)$.
Then the claim follows from an obvious identity
\[
\chi(a_1,\dots, a_{m-1},a_m+2)=\chi(a_1,\dots, a_{m-1},a_m).
\]
Let $s-i_m\>2$ and $a=\left\lfloor\frac{s-i_m}{2}\right\rfloor+1$. Then
\[
\alpha\left(n+f_s+f_{s+2}\right)=
\begin{cases}
(a_1,\dots, a_m,a,2) & \text{if\; $\beta(s)=\beta(i_m)$},\\[1mm]
(a_1,\dots, a_m)\times(a,2) & \text{if\; $\beta(s)\neq\beta(i_m)$}.
\end{cases}
\]
Now, for $\beta(s)=\beta(i_m)$, the claim follows from the equality $\xi(2;-1)=0$
and formula \eqref{eqRec}.
For $\beta(s)\neq\beta(i_m)$, it follows from the equality $\chi(a,2)=1$.
\ep

\blr\lab{hr}
Let
\[
h(i):=\left|\;\left\{n\mid 0\<n\<f_i-1\;\;\text{and}\;\;\chi(n)=0\right\}\;\right|.
\]
Then $h(0)=h(1)=h(2)=h(3)=0,h(4)=1$, and
\[
h(i)=f_{i-5}+1+h(i-1)+2h(i-4)\qquad\text{if\; $i\>5$}.
\]
\elr

\bp
For $i<5$, the claim is trivial. Let $i\>5$.

The quantity of
$n\in[f_{i-1}-1,f_i-1]$ with $\chi(n)=0$ is equal to $h(i)-h(i-1)$.
By Lemma \ref{14} we have $\chi(n)=0$ if $n\in[f_{i-4}+f_{i-1}-1,f_{i-3}+f_{i-1}]$.
The quantity of such numbers $n$ is equal to $f_{i-5}+1$.

Lemma \ref{12} implies that
\[
\left|\;\{n\in[f_{i-1}-1,f_{i-4}+f_{i-1}-2]\mid\chi(n)=0\}\;\right|=
\left|\;\{n\in[f_{i-3}+f_{i-1},f_i-1]\mid \chi(n)=0\}\;\right|. \]
Let
$n\in[f_{i-3}+f_{i-1},f_i-1]$. Clearly,
$n=n^\prime+f_{i-3}+f_{i-1}$, where $n^\prime\in[0,f_{i-4}-1]$.
Lemma \ref{15} shows that $\chi(n)=\chi(n_1)$. Thus,
\[
\left|\;\left\{n\in[f_{i-1}-1,f_{i-4}+f_{i-1}-2]\cup[f_{i-3}+f_{i-1},f_i-1]\mid\chi(n)=0\right\}\;\right|=2h(i-4).
\]
Collecting all together we obtain the asserted recurrence for $h(i)$.
\ep

\bp[Proof of Theorem \ref{th41}]
The sequence $E(a)$ does not decrease as $a$ grows, and
its subsequence $E(f_i-1)=h(i)$ is increasing.
Therefore, it suffices to check that
$\mathop{\lim}\limits_{i\to\infty}\frac{h(i)}{f_i-1}=1$.

Let $H(t):=\mathop{\sum}\limits_{i=0}^\infty h(i)t^i$.
From Lemma \ref{hr} it follows that
\[
H(t)=\frac{1}{1-t-t^2}+\frac{1}{2(t-1)}-\frac{1}{14(t+1)}+
\frac{8t^2-2t+3}{7(2t^3-2t^2+2t-1)}\;.
\]
In a standard way (see \cite{GKP}, Sec.7.3) this implies that
\[
h(i)=f_i-1+\frac{1}{7}\left(3+\beta(i)+a_1t^i_1+a_2t^i_2+a_3t^i_3\right),
\]
where $t_1,t_2,t_3$ are the roots of the equation $t^3-2t^2+2t-2=0$
and $a_1,a_2,a_3$ are some constants. Let $t_1$ be the real root while
$t_2$ and $t_3$ be the complex conjugate roots.
Then
\[
t_1=\frac{u^2+2u-2}{3u}\qquad\text{and}\qquad
t_{2,3}=\frac{2+4u-u^2}{6u}\pm\frac{2+u^2}{6u}\sqrt{-3}\,,
\qquad\mbox{where}\quad u=\sqrt[3]{17+3\sqrt{33}}\,.
\]
It is easy to verify that $|t_{1,2,3}|<\phi\approx 1.62$,
($t_1\approx 1.54,|t_{2,3}|\approx 1.13$).

Now, Theorem \ref{th41} follows from the asymptotic
$f_i\underset{i\to\infty}\sim\frac{\phi^{i+1}}{\sqrt{5}}$, see
\cite{GKP}, Sec.6.6.
\ep

\brm\lab{RemXi}
Lemma \ref{14} shows that there exist arbitrary long
intervals of consequent numbers $n$ with $\chi(n)=0$. Contrariwise,
one can show that the quantity of the consequent numbers $n$ with
$\chi(n)\neq 0$ does not exceed $4$. Moreover, only the following
sequences of $\chi(n)\neq 0$ are realized:
\[
\{1\},\{-1\},\{1,-1\},\{-1,1\},\{1,1,-1\},\{-1,-1,1\},
\{1,-1,-1\},\{-1,1,1\},\{1,-1,-1,1\},\{-1,1,1,-1\}.
\]
\erm

\section{\bf On minimal solutions of the equation $F(n)=k$ as $k$ varies}\lab{Fnk}

In this section I
analyze the problem of finding the minimal solution of the equation
$F(n)=k$ for $k>1$ from the point of view developed in
Section \ref{Gen}. The main result of the section is Theorem \ref{Fprime}.

\bdr
Denote by $m_F(k)$ the minimal solution of
the equation $F(n)=k$. We say that the number $k>1$ is
\emph{$F$-primitive} if $m_F(k)$ is $f$--simple number. In particular,
for any $F$-primitive number $k$ where is a uniquely defined
number $m(k)$ such that $\Pi(m_F(k))=\frac{m(k)}{k}$.
\edr

Let us give some examples of the $F$-primitive numbers.
For any prime $p$, any generating solution of the equation
$F(n)=p$ is $f$--simple by Proposition \ref{pr33}. Therefore, all
primes are $F$-primitive. But many composite numbers are
$F$-primitive as well. For example, Corollary \ref{Up_Est0} implies
that $m_F(f_i)=f_i^2-1$.
The number $f_i^2-1$ is $f$-simple since
\[
\Pi(f_a^2-1)=
\bcs
\pi\left\langle 2,\left\lfloor\frac{i-1}{2}\right\rfloor*3\right\rangle=\frac{f_{i-1}}{f_i} &\text{if $i\equiv 0\mod 2$},\\[1mm]
\pi\left\langle \left\lfloor\frac{i-1}{2}\right\rfloor*3\right\rangle=\frac{f_{i-2}}{f_i} &\text{if $i\equiv 1\mod 2$}.
\ecs
\]
Therefore, any Fibonacci number $f_i$ for $i\geqslant 2$ is $F$-primitive.
The first $F$-primitive numbers $k>1$ and the corresponding minimal values $m(k)$ and $m_F(k)$ are as follows:

\begin{center}
\begin{tabular}{|c||c|c|c|c|c|c|c|c|c|c|c|c|c|c|c|c|c|c|c|c|c|c|}\hline
$k$      & 2 & 3 & 5  & 7  & 8  & 11  & 13  & 17  & 18  & 19  & 21  & 23   & 27   & 29   & 31   & 34   & 37\\
\hline
$m(k)$   & 1 & 1 &  3 &  4 &  3 &   4 &   8 &  11 &   7 &  11 &   8 &   14 &   11 &   17 &   21 &   13 & 14\\
\hline
$m_F(k)$ & 3 & 8 & 24 & 58 & 63 & 152 & 168 & 406 & 401 & 435 & 440 & 1011 & 1066 & 1050 & 1160 & 1155 & 2647\\
\hline
\end{tabular}
\end{center}
\medskip

\brm
Conjecturally, the numbers $f_i+f_{i+2}$, called the \textit{Lucas numbers},
are $F$-primitive if $i>1$. Moreover,
$\Pi(m_F(f_i+f_{i+2}))=\Pi(m_F(f_i))\oplus\Pi(m_F(f_{i+2}))$,
where $\oplus$ denotes the \textit{Farey sum} (see \cite{HW}).
Then
\[
m_F(f_i+f_{i+2})=f_i^2+f_{2i+5}-1.
\]
Most likely, the numbers $f_i+f_{i+3}+f_{i+5}$ for $i\>1$ are also $F$-primitive.
\erm
Formula \eqref{eps} and the definition of $\Pi$ imply the \quat{monotonicity} of the multiplication $\times$:

\blr\lab{lm83}
For any $n,n_1,n_2\in\mathbb{G}$ such that $n_1<n_2$ we have $n_1\times n<n_2\times n$.
\elr

\btr\lab{Fprime}
For any natural $k\>2$ let
\[
\Pi(m_F(k))=\frac{l_1}{k_1}\times\frac{l_2}{k_2}\times\dots\times\frac{l_r}{k_r}\in\Gamma_+.
\]
Then the numbers $k_1,k_2,\dots,k_r$ are $F$-primitive and
$l_i=m(k_i)$ for any $i\in[1,r]$. Thus
\[
m_F(k)=m_F(k_1)\times m_F(k_2)\times\dots\times m_F(k_r).
\]
\etr
\bp
The claim follows from Lemma \ref{lm83} by induction on $r$.
\ep

For example,
\begin{align*}
m_F(6)&=m_F(2)\times m_F(3)=37,&
m_F(54)&=m_F(3)\times m_F(18)=4456,&\\[1mm]
m_F(57)&=m_F(19)\times m_F(3)=4616,&
m_F(96)&=m_F(2)\times m_F(2)\times m_F(8)\times m_F(3)=12093.&
\end{align*}
Theorem \ref{Up_Est_r} easily implies that\;
$m_F(2f_i)=m_F(2)\times m_F(f_i)=3\times(f_i^2-1)=f_if_{i+3}-2\beta(i)$.

The following questions naturally arise:
\en
\item
How to describe the set of $F$-primitive numbers?
\item
How to find $m(k)$ for an $F$-primitive number $k$ (for example, when $k$ is prime)?
\item
When for a sequence of $F$-primitive numbers $k_1,k_2,\dots,k_r$ we have
\[
m_F(k_1k_2\cdots k_r)=m_F(k_1)\times m_F(k_2)\times\dots\times m_F(k_r)?
\]
\ene

I have plotted the graph of linear approximation
of the frequencies distribution (histogram) of rational numbers $\frac{m(k)}{k}$,
where $k$ runs over the set of $F$-primitive numbers from the segment $[1,5000]$.
This segment contains 1764 such numbers, 669 of them are primes.
To build the plot, the sequence of minimal solutions of the equation $F(n)=k$ from \cite{SL}, A013583,
was quite useful.

\begin{figure}[h!]
  \centering
  \includegraphics[width=14cm,keepaspectratio]{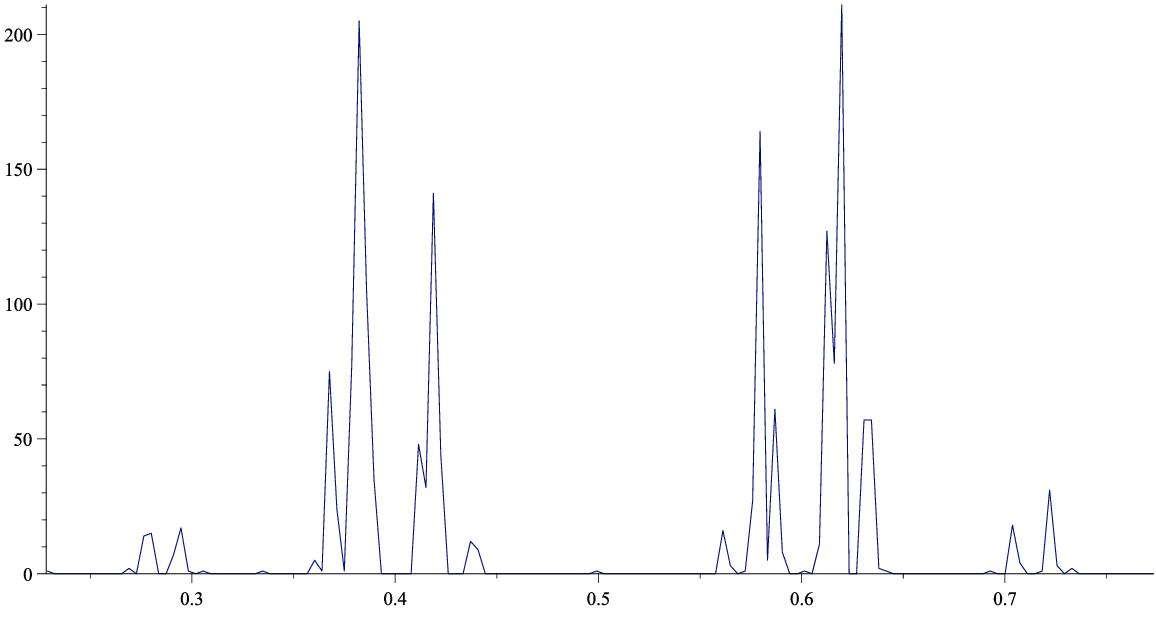}
  \caption{\small\sf The frequency distribution of rationals $\Pi(m_F(k))=\frac{m(k)}{k}$, where $k$ runs
  over the set of $F$-primitive numbers in interval $[1,5000]$.}
  \label{Fig1}
\end{figure}

The graph shows that set $\mathbb{P}$ of rational numbers $\frac{m(k)}{k}$,
where $k$ runs over the set of all $F$-primitive numbers, has a complicated structure.
Probably, this set has infinite number of limit points.

\bhr\lab{conj}
The set $\mathbb{P}$ is invariant with respect
to the reflection of the interval $[0,1]$ with respect to point $\frac{1}{2}$.
The numbers $\phi^{-2}$ and $\phi^{-1}$,
where $\phi=\frac{1}{2}\left(1+\sqrt{5}\right)$,
are the limit points of $\mathbb{P}$.
\ehr

\section{\bf On the graph of the function $F$}\lab{Fgr}

Many of the claims in this article are the result of attempting
to explain the behavior of the graph of $F$.
Let me make several remarks on this. The graph on the
interval $[1,f_{19}-1]$ is presented below:

\begin{figure}[h!]
  \centering
  \includegraphics[width=16cm,keepaspectratio]{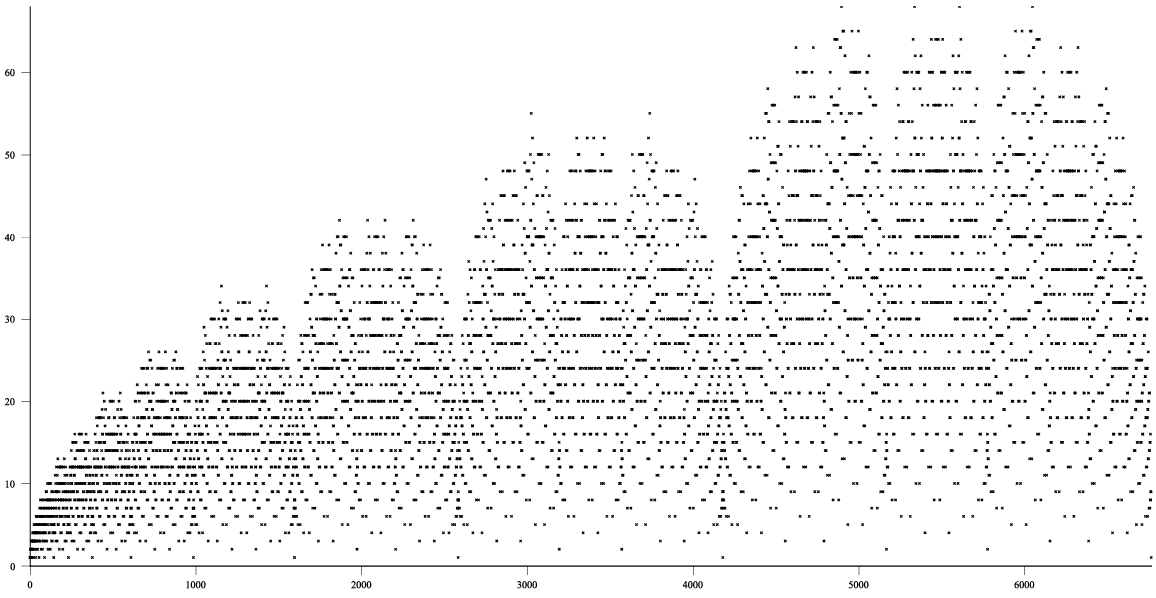}
  \caption{\small\sf The graph of the function $F(n)$, where $n\in[1,f_{19}-1]$.}
  \label{figF}
\end{figure}

Theorem \ref{th_F} implies that $F(n)=1$ whenever
$n=f_i-1$, where $i\geqslant 1$. The graph of $F$ consists of a sequence
of increasing \quat{waves}: the $i$th wave $W_i$ for $i\geqslant 2$ is the set
of points in the plane with coordinates $(n, F(n))$, where
$f_i\leqslant n<f_{i+1}$; by definition, $W_0=(0,1),W_1=(1,1)$.
For $f_a\leqslant n<f_{a+1}$, where $a\leqslant i-2$, formula \eqref{PhiTi} for
$t=1$ shows that
\beq\lab{n+}
F(n+f_i)=
\bcs
\frac{i-a+1}{2}\cdot F(n)&\text{if $a\not\equiv i\mod 2$},\\
\frac{i-a+2}{2}\cdot F(n)-F(n-f_a)&\text{if $a\equiv i\mod 2$}.
\ecs
\eeq
Let $S_i(n, F(n)):=(n+f_i,F(n+f_i))$. It is easy to see
that
\beq\lab{Wi}
W_i = S_i(W_0\cup W_1\cup\dots\cup W_{i-2}),\qquad\text{where $i\geqslant 2$}.
\eeq
This allows us to build the sets $W_i$ inductively, starting
from the set $W_0\cup W_1$. Formulas \eqref{Wi} and \eqref{n+} show
that the $i$th wave, where $i\geqslant 2$, is “similar” to the union of
the $i-2$ previous waves. These formulas can also be used
to explain some visible symmetries inside each wave.

Lemma \ref{12} shows that the function $F$ is invariant with
respect to the refection $\rho_i$ of the interval $[f_i-1,f_{i+1}-1]$
with respect to its center.

The next conjecture can provide a much more
precise upper bound for $F(n)$ when $n\in[f_i-1,f_{i+1}-1]$
as compared with the estimate in Theorem \ref{Up_Est_r}:

\bhr\lab{conjc}
Let $C_i$ be the convex hull of the
set of points $(n, F(n))$, where $n\in[f_i-1, f_{i+1}-1]$ and
$i>7$. Set
\begin{gather*}
b_1(i,q):=f_i-1+f_q^2,\qquad b_2(i,q):=f_i-1+f_qf_{q+1}+2\cdot(-1)^{q+1}, \\
B_r(i):=\left\{b_r(i,q),\rho_i(b_r(i,q))\mid 1\leqslant q\leqslant\left\lfloor\frac{i}{2}\right\rfloor-2\right\},
\end{gather*}
where $r=1,2$ and $\rho_i(n)=f_{i+2}-2-n$.

Then $C_i$ coincides with the convex hull of the set
of points $(c,F(c))$, where $c\in B_1(i)\cup\{f_i-1,f_{i+1}-1\}$
if $\beta(i)=1$, and $c\in B_1(i)\cup B_2(i)\cup\{f_i-1, f_{i+1}-1\}$ if
$\beta(i)=0$.
\ehr

\section{\bf Variations on the theme of function $\Psi(k)$}\lab{Add}

The material of this section concerns the function
$\Psi(k)=|\mathbb{G}(k)|$
as well as several additional functions similar to it.
In particular, from the formulas we will obtain
a recurrent expression for $\Psi(k)$
and a formula for the Dirichlet generating function of the sequence $\Psi(k)$ will follow.

Let $k=p_1^{n_1}p_2^{n_2}\cdots p_r^{n_r}$ be the prime
decomposition of $k\in\mathbb{N}$ and $\deg(k):=n_1+n_2+{\dots}+n_r$.
Denote by $\G_+(k,m)$ the set of words
$\frac{l_1}{k_1}\times\frac{l_2}{k_2}\times{\dots}\times\frac{l_m}{k_m}\in\G_+(k)$.
Then
\[
\G_+=\coprod_{k=2}^\infty\;\coprod_{m=1}^{\deg(k)}\G_+(k,m)
\]
is a disjoint union of nonempty subsets. Set
\[
\psi(k,m):= \bcs
\left|\G_+(k,m)\right|&\text{if $m\<\deg(k)$},\\
0&\text{if $m>\deg(k)$} \ecs \qquad\text{and}\qquad
\Psi(k;t):=\sum_{m=1}^\infty\psi(k,m)t^m.
\]

\brm
In Section \ref{Add}, by abuse of notation I denote the function
of two arguments by the same symbol $\Psi$ used to denote
the function of one argument $\Psi(k):=\Psi(k;1)$.
\erm

\btr
The function $\Psi(k;t)$ is expressed by the recursive formula
\[
\Psi(k;t)=t\sum_{a>1,\;a|k}\Psi\left(\frac{k}{a}\,;\,t\right)\,\varphi(a),
\]
where $\Psi(1;t)=1$.
\etr

\bp
Let $a>1$ and
$\G^{(a)}_+(k,m):=\left\{\frac{l_1}{k_1}\times\frac{l_2}{k_2}\times{\dots}\times\frac{l_m}{k_m}\in\G_+(k)\mid k_m=a\right\}$.
Then $\G^{(a)}_+(k,m)\neq\emptyset$ if and only if $a|k$, and
\beq\lab{gam}
\left|\G^{(a)}_+(k,m)\right|=
\sum_{1\<b<a,\;\gcd(a,b)=1}\left|\G_+\left(\frac{k}{a}\,,\,m-1\right)\times\frac{b}{a}\right|=
\left|\G_+\left(\frac{k}{a}\,,\,m-1\right)\right|\cdot\varphi(a).
\eeq
Since $\G^{(a_1)}_+(k,m)\cap\G^{(a_2)}_+(k,m)=\emptyset$ for $a_1\neq a_2$, it follows that
\beq\lab{eqq}
\left|\G_+(k,m)\right|=\sum_{a>1,\;a|k}\left|\G^{(a)}_+(k,m)\right|.
\eeq
To finish the proof we substitute  expression \eqref{gam} in \eqref{eqq}
and apply the definition of $\Psi(k;t)$.
\ep

\btr\lab{dgft}
We have
\beq\lab{dir_psi}
1+\sum_{k=2}^\infty\,\frac{\Psi(k;t)}{k^s}=\left(1+t-t\;\frac{\zeta(s-1)}{\zeta(s)}\right)^{-1}.
\eeq
\etr

\bp(Cf.\cite{Hill}, p.138.)
Let
\[
L(s,\varphi)=1+\sum_{k=2}^\infty\,\frac{\varphi(k)}{k^s}\;.
\]
For any $m\>1$, we obviously have
\[
\sum_{k=2}^\infty\frac{\psi(k,m)}{k^s}=\left(L(s,\varphi)-1\right)^m.
\]
Therefore,
\[
1+\sum_{k=2}^\infty\,\frac{\Psi(k;t)}{k^s}=1+\sum_{k=2}^\infty\sum_{m=1}^\infty\,\frac{\psi(k,m)t^m}{k^s}=
1+\sum_{m=1}^\infty t^m\left(L(s,\varphi-1)\right)^m=\left(1+t-tL(s,\varphi)\right)^{-1}.
\]
Now, the claim required follows from the classical formula
$L(s,\varphi)=\frac{\zeta(s-1)}{\zeta(s)}\;$ (see \cite{HW}, p.250).
\ep

Formula \eqref{dir_psi} for $t=1$ turns into formula \eqref{zet}.
For $t=-1$, formula \eqref{dir_psi} gives, after straightforward calculations, the expression
\[
\Psi(k;-1)=\mu(k)\,\varphi(k),
\]
where $\mu$ is the \emph{M\"obius function} (see \cite{HW}, p.234).

The sequence $\Psi(k)$ rapidly grows with the growth of the quantity of prime factors of $k$.
Let $p$ be a prime, let $p_1,p_2,\dots, p_r$ be different primes.
One can show that
\[
\Psi(p^n)=(p-1)(2p-1)^{n-1}\qquad\text{and}\qquad
\Psi(p_1p_2\cdots p_r)=B(r)\;(p_1-1)(p_2-1)\cdots(p_r-1),
\]
where $B(0)=1$ and
\[
B(r+1)=\sum_{i=0}^r\binom{r+1}{i}B(i).
\]
The numbers $B(r)\in\{1,3,13,75,541,4683,47293,545835,\cdots\}$, where $r\>1$,
are known as \textit{ordered Bell numbers} (see \cite{Wilf}, p.175).
The number $B(r)$ is the quantity of ways to represent the set
$\{1,2,\dots, r\}$ as an ordered union of non-intersecting subsets.

\bdr
The words
$\frac{l_1}{k_1}\times\frac{l_2}{k_2}\times{\dots}\times\frac{l_m}{k_m}$ and
$\frac{l^\prime_1}{k^\prime_1}\times\frac{l^\prime_2}{k^\prime_2}\times{\dots}\times\frac{l^\prime_m}{k^\prime_m}$
from the set $\G_+(k,m)$ will be referred to as
\vspace{1mm}

\textit{$\omega$-equivalent} if $k_1=k^\prime_1,k_2=k^\prime_2,\dots,k_m=k^\prime_m$;
\vspace{1mm}

\textit{$t$-equivalent} if
$\left\{\frac{l_1}{k_1},\frac{l_2}{k_2},\dots,\frac{l_m}{k_m}\right\}=
\left\{\frac{l^\prime_1}{k^\prime_1},\frac{l^\prime_2}{k^\prime_2},\dots,\frac{l^\prime_m}{k^\prime_m}\right\}$ as sets;
\vspace{1mm}

\textit{$s$-equivalent} if
$\left\{k_1,k_2,\dots,k_m\right\}=\left\{k^\prime_1,k^\prime_2,\dots,k^\prime_m\right\}$ as sets;
\vspace{1mm}

\textit{$*$-equivalent} if they are $t$-equivalent, and if
$k_i=k^\prime_i$, then either $l_i=l^\prime_i$, or $l_i\cdot
l^\prime_i\equiv 1\mod k_i$.
\edr

Let us extend each of these relations to $\G_+$ component-wise. It
is easy to see that, for any of the relations introduced, the elements
of $\G_+$ are equivalent if and only if they belong to an orbit of
the action of their corresponding subgroups in the automorphism group of $\G_+$.

Each \textit{$\omega$-equivalence} class from the set $\G_+(k,m)$
contains a unique word
$\frac{1}{k_1}\times\frac{1}{k_2}\times{\dots}\times\frac{1}{k_m}$,
where $k_1 k_2\cdots k_m=k$, and the order of factors
is essential. Any such decomposition is called an \textit{ordered
multiplicative partition} of $k$ of length $m$. Let
$\widehat{\psi}(k,m)$ be the quantity of such partitions. Define
\[
\widehat{\Psi}(k;t):=\sum_{m=1}^\infty\widehat{\psi}(k,m)t^m.
\]
Arguing exactly as in the proof of Theorem \ref{dgft},
where $\varphi$ is replaced with $1$, we obtain the formula
\[
1+\sum_{k=2}^\infty\,\frac{\widehat{\Psi}(k;t)}{k^s}=\left(1+t-t\zeta(s)\right)^{-1},
\]
which is well known at least for $t=1$, see \cite{Hill}.

Each \textit{$t$-equivalence} class from the set $\G_+(k,m)$ contains
a unique word
$\frac{l_1}{k_1}\times\frac{l_2}{k_2}\times{\dots}\times\frac{l_m}{k_m}$,
such that $1<k_1\<k_2\<{\dots}\<k_m$, and if $k_i=k_{i+1}$, then
$l_i\<l_{i+1}$ for $1\<i<m$. The set of these classes has a natural
structure of a free \textit{commutative monoid} $\G_+^c$, generated
by the set $\mathbb{Q}_{(0,1)}$.

Denote the set of its elements with $k_1 k_2\cdots k_r=k$ by $\G_+^c(k)$ and define
$\Psi_c(k):=\left|\G_+^c(k)\right|$. One can think about $\Psi_c$ as about
a \quat{commutative} analog of $\Psi$. To obtain a generating
function for the sequence $\Psi_c(k)$, let us first make a general
observation.

For any function $f:\mathbb{N}\to\mathbb{N}$, consider the
decomposition
\[
\prod_{k=2}^\infty\left(1-\frac{1}{k^s}\right)^{-f(k)}=
1+\sum_{k=2}^\infty\frac{\widetilde{f}(k)}{k^s}\;.
\]
One can interpret the number $\widetilde{f}(k)$ as the number of
ways to write $k$ as the product $k=k_1 k_2\cdots k_r$,
where $1<k_1\<k_2\<{\dots}\<k_r$, each factor $k_i$ has $f(k_i)$
different colors, and differently colored products are considered
different. One can show that, for a prime $p$, \beq\lab{pn}
\widetilde{f}(p^a)=\sum_{a_1+2a_2+{\dots}+ sa_s+{\dots}=a}\;
\prod_{i=1}^\infty\binom{f(p^i)+a_i-1}{a_i}, \eeq where $a_1,a_2,\dots,
a_s,{\dots}$ are the non-negative integers.

Let $f(k)=\varphi(k)$ be the Euler totient function. Then, obviously,
\[
\prod_{k=2}^\infty\left(1-\frac{1}{k^s}\right)^{-\varphi(k)}=1+\sum_{k=2}^\infty\frac{\Psi_c(k)}{k^s}\;.
\]
One can show that, for different primes $p_1,p_2,\dots, p_r$, we have
\[
\Psi_c(p_1p_2\cdots p_r)=B_c(r)\,(p_1-1)(p_2-1)\cdots(p_r-1),
\]
where $B_c(0)=1$ and
\[
B_c(r+1)=\sum_{i=0}^r\binom{r}{i}B_c(i).
\]
The numbers $B_c(r)\in\{1,2,5,15,52,203,877,4140,{\dots}\}$, where $r\>1$,
are known as \textit{Bell numbers} (see \cite{Wilf}, p.20).
The number $B_c(r)$ is the quantity of ways to represent the set
$\{1,2,\dots, r\}$ as a union of non-intersecting subsets.

I was unable to find a general expression for $\Psi_c(p^n)$.
Here are the first several ones:
\begin{eqnarray*}
\Psi_c(p)\;\,   &=& p-1,\\
\Psi_c(p^2) &=& \frac{p(p-1)}{2!}\cdot 3,\\
\Psi_c(p^3) &=& \frac{p(p-1)}{3!}\cdot (13p-5),\\
\Psi_c(p^4) &=& \frac{p(p-1)}{4!}\cdot (73p^2-45p+14),\\
\Psi_c(p^5) &=& \frac{p(p-1)}{5!}\cdot (501p^3-414p^2+111p-54),\\
\Psi_c(p^6) &=& \frac{p(p-1)}{6!}\cdot (4051p^4-4130p^3+1445p^2-190p+264).
\end{eqnarray*}

For the case of \textit{$s$-equivalence}, we set
$f(k)=1$ for any $k\>2$. Then we obtain the known formula
\[
\prod_{k=2}^\infty\left(1-\frac{1}{k^s}\right)^{-1}=1+\sum_{k=2}^\infty\frac{\widehat{\Psi}_c(k)}{k^s}\;,
\]
where $\widehat{\Psi}_c(k)$ is the quantity of the decompositions
$k=k_1 k_2\cdots k_m$, where $1<k_1\<k_2\<{\dots}\<k_m$.
Such a decomposition is called an \textit{unordered multiplicative partition} of $k$.
For example, $\widehat{\Psi}_c(p^n)$ for a prime $p$ is equal to the
quantity of all integer partitions of $n$.
This evident claim follows from formula \eqref{pn} as well.

Finally, the generating function for the number
of \textit{$*$-equivalent} classes, has the form
\[
\prod_{k=2}^\infty\left(1-\frac{1}{k^s}\right)^{-\varphi_*(k)}=
1+\sum_{k=2}^\infty\frac{\Psi_*(k)}{k^s}\;,
\]
where $\varphi_*(k)$ is defined as follows. Let
$(\mathbb{Z}/k\mathbb{Z})^*$ denote the multiplicative group of the
ring $\mathbb{Z}/k\mathbb{Z}$. Then $\varphi_*(k)$ is the quantity
of the invariant subsets under the inversion operation in the group
$(\mathbb{Z}/k\mathbb{Z})^*$.

The function $\varphi_*$ is multiplicative. That is,
$\varphi_*(k_1 k_2)=\varphi_*(k_1)\,\varphi_*(k_2)$ for relatively prime $k_1$ and $k_2$.
For any prime $p$ and any $a\in\mathbb{N}$,
the structure of the group $(\mathbb{Z}/p^a\mathbb{Z})^*$ is well known (see \cite{IR}).
It is not difficult to check the formula
\[
\varphi_*(p^a)=
\begin{cases}
1 & \text{if $p=2$ and $a=1$},\\
2 & \text{if $p=2$ and $a=2$},\\
2^{a-2}+2 & \text{if $p=2$ and $a>2$},\\
\frac{1}{2}\, p^{a-1}\, (p-1)+1 & \text{if $p>2$}.
\end{cases}
\]
This allows us to easily compute $\varphi_*(n)$ for any natural $n$.

Calculations show that the first $20$ values of the functions
$\Psi(k),\Psi_c(k),\Psi_*(k)$ are as follows (see \cite{SL},
\,A006874,\,A007896,\,A007898):
\medskip

\begin{center}
\begin{tabular}{|c||c|c|c|c|c|c|c|c|c|c|c|c|c|c|c|c|c|c|c|c|c|c|c|c|c|}\hline
$k$ & 1 & 2 & 3 & 4 & 5 & 6 & 7 & 8 & 9 & 10 & 11 & 12 & 13 & 14 & 15 & 16 & 17 & 18 & 19 & 20\\
\hline
\hline
$\Psi(k)$ & 1 & 1 & 2 & 3 & 4 & 6 & 6 & 9 & 10 & 12 & 10 & 22 & 12 & 18 & 24 & 27 & 16 & 38 & 18 & 44\\
\hline
$\Psi_c(k)$ & 1 & 1 & 2 & 3 & 4 & 4 & 6 & 7 & 9 & 8 & 10 & 12 & 12 & 12 & 16 & 18 & 16 & 19 & 18 & 24\\
\hline
$\Psi_*(k)$ & 1 & 1 & 2 & 3 & 3 & 4 & 4 & 7 & 7 & 6 & 6 & 12 & 7 & 8 & 12 & 16 & 9 & 15 & 10 & 18\\
\hline
\end{tabular}
\end{center}


\end{document}